
\def\title#1{{\titlefont\noindent #1\bigskip}}

\def\author#1{{\largefont\noindent #1}\medskip}

\def\beginlinemode{\endmode
 \begingroup\obeylines\def\endmode{\par\endgroup}}
\let\endmode=\par

\newbox\theaddress
\def\address{\smallskip\beginlinemode\parindent 0in\getaddress}
{\obeylines
\gdef\getaddress #1 
 #2
 {#1\gdef\addressee{#2}%
   \global\setbox\theaddress=\vbox\bgroup\raggedright%
    \everypar{\hangindent2em}#2
   \def\endaddress{\egroup\endgroup \copy\theaddress \medskip}}}

\def\thanks#1{\footnote{}{\eightpoint #1}}

\long\def\Abstract#1{{\eightpoint\narrower\vskip\baselineskip\noindent
#1\smallskip}}

\def\skipfirstword#1 {}

\def\ir#1{\csname #1\endcsname}

\newdimen\currentht
\newbox\droppedletter
\newdimen\droppedletterwdth
\newdimen\drophtinpts
\newdimen\dropindent

\def\irrnSection#1#2{\edef\tttempcs{\ir{#2}}
\currentht-\pagetotal\advance\currentht by-\ht\footins
\advance\currentht by\vsize
\ifdim\currentht<1.5in\par\vfill\eject\else\vbox to.25in{\vfil}\fi
{\largefont\noindent{}\expandafter\skipfirstword\tttempcs. #1}
\vskip\baselineskip }

\def\irSubsection#1#2{\edef\tttempcs{\ir{#2}}
\vskip\baselineskip\penalty-3000
{\bf\noindent \expandafter\skipfirstword\tttempcs. #1}
\vskip\baselineskip}

\def\irSubsubsection#1#2{\edef\tttempcs{\ir{#2}}
\vskip\baselineskip\penalty-3000
{\bf\noindent \expandafter\skipfirstword\tttempcs. #1}
\vskip\baselineskip}

\def\References{\vbox to.25in{\vfil}\noindent{}{\bf References}
\vskip\baselineskip\par}

\def\baselinebreak{\par \ifdim\lastskip<\baselineskip
         \removelastskip\penalty-200\vskip\baselineskip\fi}

\long\def\prclm#1#2#3{\baselinebreak
\noindent{\bf \csname #2\endcsname}:\enspace{\sl #3\par}\baselinebreak}

\def\Prf{\noindent{\bf Proof}: }

\def\rem#1#2{\baselinebreak\noindent{\bf \csname #2\endcsname}:\enspace }

\def\qed{{\hfill$\diamondsuit$}\vskip\baselineskip}

\def\bibitem#1{\par\indent\llap{\rlap{\bf [#1]}\indent}\indent\hangindent
2\parindent\ignorespaces}

\long\def\eatit#1{}

\def\leftheadlinetext{}
\def\rightheadlinetext{}

\def\leftheadline{{\eightrm\folio\hfil \leftheadlinetext\hfil}}
\def\rightheadline{{\eightrm\hfil\rightheadlinetext\hfil\folio}}

\headline={\ifnum\pageno=1\hfil\else
\ifodd\pageno\rightheadline\else\leftheadline\fi\fi}

\def\tenpoint{\def\rm{\fam0\tenrm}
\textfont0=\tenrm \scriptfont0=\sevenrm \scriptscriptfont0=\fiverm
\textfont1=\teni \scriptfont1=\seveni \scriptscriptfont1=\fivei
\def\mit{\fam1} \def\oldstyle{\fam1\teni}
\textfont2=\tensy \scriptfont2=\sevensy \scriptscriptfont2=\fivesy
\def\cal{\fam2}
\textfont3=\tenex \scriptfont3=\tenex \scriptscriptfont3=\tenex
\def\it{\fam\itfam\tenit} 
\textfont\itfam=\tenit
\def\sl{\fam\slfam\tensl} 
\textfont\slfam=\tensl
\def\bf{\fam\bffam\tenbf} 
\textfont\bffam=\tenbf \scriptfont\bffam=\sevenbf
\scriptscriptfont\bffam=\fivebf
\def\tt{\fam\ttfam\tentt} 
\textfont\ttfam=\tentt
\normalbaselineskip=12pt
\setbox\strutbox=\hbox{\vrule height8.5pt depth3.5pt  width0pt}%
\normalbaselines\rm}

\def\eightpoint{\def\rm{\fam0\eightrm}%
\textfont0=\eightrm \scriptfont0=\sixrm \scriptscriptfont0=\fiverm
\textfont1=\eighti \scriptfont1=\sixi \scriptscriptfont1=\fivei
\def\mit{\fam1} \def\oldstyle{\fam1\eighti}%
\textfont2=\eightsy \scriptfont2=\sixsy \scriptscriptfont2=\fivesy
\def\cal{\fam2}%
\textfont3=\tenex \scriptfont3=\tenex \scriptscriptfont3=\tenex
\def\it{\fam\itfam\eightit} 
\textfont\itfam=\eightit
\def\sl{\fam\slfam\eightsl} 
\textfont\slfam=\eightsl
\def\bf{\fam\bffam\eightbf} 
\textfont\bffam=\eightbf \scriptfont\bffam=\sixbf
\scriptscriptfont\bffam=\fivebf
\def\tt{\fam\ttfam\eighttt} 
\textfont\ttfam=\eighttt
\normalbaselineskip=9pt%
\setbox\strutbox=\hbox{\vrule height7pt depth2pt  width0pt}%
\normalbaselines\rm}

\def\largefont{\def\rm{\fam0\largerm}
\textfont0=\largerm \scriptfont0=\largescriptrm \scriptscriptfont0=\tenrm
\textfont1=\largei \scriptfont1=\largescripti \scriptscriptfont1=\teni
\def\mit{\fam1} \def\oldstyle{\fam1\teni}
\textfont2=\largesy 
\def\cal{\fam2}
\def\it{\fam\itfam\largeit} 
\textfont\itfam=\largeit
\def\sl{\fam\slfam\largesl} 
\textfont\slfam=\largesl
\def\bf{\fam\bffam\largebf} 
\textfont\bffam=\largebf 
\def\tt{\fam\ttfam\largett} 
\textfont\ttfam=\largett
\normalbaselineskip=17.28pt
\setbox\strutbox=\hbox{\vrule height12.25pt depth5pt  width0pt}%
\normalbaselines\rm}

\def\titlefont{\def\rm{\fam0\titlerm}
\textfont0=\titlerm \scriptfont0=\largescriptrm \scriptscriptfont0=\tenrm
\textfont1=\titlei \scriptfont1=\largescripti \scriptscriptfont1=\teni
\def\mit{\fam1} \def\oldstyle{\fam1\teni}
\textfont2=\titlesy 
\def\cal{\fam2}
\def\it{\fam\itfam\titleit} 
\textfont\itfam=\titleit
\def\sl{\fam\slfam\titlesl} 
\textfont\slfam=\titlesl
\def\bf{\fam\bffam\titlebf} 
\textfont\bffam=\titlebf 
\def\tt{\fam\ttfam\titlett} 
\textfont\ttfam=\titlett
\normalbaselineskip=24.8832pt
\setbox\strutbox=\hbox{\vrule height12.25pt depth5pt  width0pt}%
\normalbaselines\rm}

\nopagenumbers

\font\eightrm=cmr8
\font\eighti=cmmi8
\font\eightsy=cmsy8
\font\eightbf=cmbx8
\font\eighttt=cmtt8
\font\eightit=cmti8
\font\eightsl=cmsl8
\font\sixrm=cmr6
\font\sixi=cmmi6
\font\sixsy=cmsy6
\font\sixbf=cmbx6

\font\largerm=cmr12 at 17.28pt
\font\largei=cmmi12 at 17.28pt
\font\largescriptrm=cmr12 at 14.4pt
\font\largescripti=cmmi12 at 14.4pt
\font\largesy=cmsy10 at 17.28pt
\font\largebf=cmbx12 at 17.28pt
\font\largett=cmtt12 at 17.28pt
\font\largeit=cmti12 at 17.28pt
\font\largesl=cmsl12 at 17.28pt

\font\titlerm=cmr12 at 24.8832pt
\font\titlei=cmmi12 at 24.8832pt
\font\titlesy=cmsy10 at 24.8832pt
\font\titlebf=cmbx12 at 24.8832pt
\font\titlett=cmtt12 at 24.8832pt
\font\titleit=cmti12 at 24.8832pt
\font\titlesl=cmsl12 at 24.8832pt

\tenpoint



\hsize 6.5in
\vsize 9in
\tolerance 3000
\hbadness 3000

\def\trans{6}
\def\duke{7}
\def\vanc{8}
\def\ravello{9}
\def\mtnwest{10}
\def\ars{11}
\def\bir{12}
\def\syzconf{13}

\def\C#1{\hbox{$\cal #1$}}
\def\s{\hbox{$\scriptstyle\cal S$}}
\def\r{\hbox{$\scriptstyle\cal R$}}

\def\C#1{\hbox{$\cal #1$}}
\def\pr#1{\hbox{{\bf P}${}^{#1}$}}

\def\smallbinom#1#2{({{\scriptstyle #1}\atop{\scriptstyle #2}})}

\newbox\mylimboxa
\newbox\mylimboxb

\newdimen\limdimht
\newdimen\limdimwd

\def\mylim#1#2#3{\setbox\mylimboxa=\hbox{#1}%
\setbox\mylimboxb=\hbox{#2}%
\limdimht\wd\mylimboxa\limdimwd\wd\mylimboxb
\ifdim\limdimht>\limdimwd\limdimwd\limdimht\fi
\setbox\mylimboxa=\hbox to\limdimwd{\hfil #1\hfil}%
\setbox\mylimboxb=\hbox to\limdimwd{\hfil #2\hfil}%
\limdimht\ht\mylimboxb\advance\limdimht by #3.0pt%
\lower\limdimht\box\mylimboxb\kern-\limdimwd\box\mylimboxa\hbox{$\,$}}

\def\leftheadline{{\eightrm\folio\hfil \leftheadlinetext\hfil}}
\def\rightheadline{{\eightrm\hfil\rightheadlinetext\hfil\folio}}

\def\leftheadlinetext{Brian Harbourne}
\def\rightheadlinetext{Resolutions of Fat Point Ideals}

\title{Free Resolutions of Fat Point Ideals on \pr2}

\author{Brian Harbourne}

\address
Department of Mathematics and Statistics
University of Nebraska-Lincoln
Lincoln, NE 68588-0323
email: bharbourne@unl.edu
WEB: http://www.math.unl.edu/$\sim$bharbour/
\smallskip
August 6, 1995\endaddress
\vskip-\baselineskip

\thanks{\vskip -6pt
\noindent This work was supported both by the National Science Foundation and
by a Spring 1994 University of Nebraska Faculty Development Leave,
and is a direct outgrowth of a visit to the Curves Seminar at Queens.
\smallskip
\noindent {\it Mathematics Subject Classification. } Primary 13D02, 14C20. 
Secondary 14J26, 13D40, 13P99.
\smallskip
\noindent {\it Key words and phrases. }  Free resolution, 
complete ideal, essentially distinct points, fat points, rational surface.\smallskip}

\Abstract{Minimal free resolutions for homogeneous ideals corresponding
to certain 0-dimensional subschemes of \pr2 
defined by sheaves of complete ideals are determined implicitly. 
All work is over an algebraically closed field 
of arbitrary characteristic.}

\irrnSection{Introduction}{introsect}
Given distinct points $p_1,\ldots,p_r$ of a smooth variety $V$
(over an algebraically closed field $k$)
and positive integers $m_i$, $Z=m_1p_1+\cdots+m_rp_r$
denotes the subscheme defined locally at each point $p_i$
by $\C I_i^{m_i}$, where $\C I_i$ is the maximal ideal in the local ring
$\C O_{V,p_i}$ at $p_i$ of the structure sheaf. More briefly,
we say $Z$ is a {\it fat point subscheme\/} of $V$. In the case that $V$
is \pr{n} for some $n$, it is of interest to study the homogeneous ideal
$I_Z$ defining $Z$ as a subscheme of \pr{n}; $I_Z$ is called an 
{\it ideal of fat points}. 

Given an ideal $I_Z$ of a fat point subscheme $Z\subset\pr{n}$, one first may want to
determine its Hilbert function $h_{I_Z}$, defined for each $d$ by 
$h_{I_Z}(d)=\hbox{dim}_k((I_Z)_d)$, where $(I_Z)_d$ is the homogeneous component
of $I_Z$ of degree $d$ (i.e., $(I_Z)_d$ is the $k$-vector space of all
homogeneous forms of degree $d$ in $I_Z$). One next may wish to study
a minimal free resolution $\cdots \to F_1\to F_0\to I_Z\to 0$ of $I_Z$,
beginning with determining $F_0$. Determining $F_0$ as a graded module
over the homogeneous coordinate ring $k[\pr n]$ of \pr n is equivalent to
finding the number $\nu_d(Z)$ of generators in each degree $d$ in 
a minimal homogeneous set of generators for $I_Z$, since
$F_0=\bigoplus_{d}R[-d]^{\nu_d}$ (where $R$ denotes $k[\pr n]$ graded by total degree,
and $R[-d]$ signifies that the degree has been shifted such that constants have degree $d$).

The Hilbert function $h_{I_Z}$ in the case that the points $p_1,\ldots,p_r\subset\pr n$ are general
has attracted attention (see [16] in general, or  [\vanc], [15], [5] and [\ravello] for $n=2$)
but much remains conjectural. Most work done on minimal free resolutions
of $I_Z$ has been restricted to the case that $Z$ is a smooth union of general points
(cf. [17]). More can be said in the case of subschemes of \pr2 involving 
small numbers of points or points in special position. 
For example, by [\ars], $h_{I_Z}$ is completely 
understood for any fat point subscheme $Z=m_1p_1+\cdots+m_rp_r\subset\pr2$ where
$p_1,\ldots,p_r$ are points (even possibly infinitely near) of a plane cubic (possibly 
reducible and nonreduced), by [2] one can determine a minimal homogeneous set of generators
for $I_Z$ in case $p_1,\ldots,p_r$ lie on a smooth plane conic, and by [\syzconf]
one can determine a minimal homogeneous set of generators
for $I_Z$ in case $Z=m(p_1+\cdots+p_r)$, where $p_1,\ldots,p_r$ are $r\le 9$ general points
of \pr2 ([\syzconf] also conjectures a result for $r>9$).

Here we will be concerned with determining minimal homogeneous sets of generators
for ideals $I_Z$ where $Z=m_1p_1+\cdots+m_rp_r\subset\pr2$ in the case that
$p_1,\ldots,p_r$ lie on a plane curve of degree at most 3. We obtain results in some special cases
in the case of a smooth curve of degree 3, but we obtain complete results
for curves of degrees 1 and 2. (Our results for points on a conic are distinguished from those 
of [2] in that we consider arbitrary conics and we allow infinitely near points.)

Since a fat point ideal $I_Z$ is perfect, one feature of working on \pr2 is that
a minimal free resolution of $I_Z$ is
of the form $0\to F_1\to F_0\to I_Z\to 0$ for some graded free modules
$F_0$ and $F_1$. Moreover, $F_0$ and $F_1$ are determined
by the Hilbert function $h_{I_Z}$ of $I_Z$
and a minimal homogeneous set of generators for $I_Z$.
In particular, we saw above that $F_0=\bigoplus_{d}R[-d]^{\nu_d}$
(where here $R=k[\pr2]$)
which means the Hilbert function $h_{F_0}$ of $F_0$ is 
$h_{F_0}(n)=\sum_{d}\nu_d\smallbinom{n+2-d}{2}$.
But the Hilbert function of
$F_1$ is determined by the Hilbert functions of $I_Z$ and $F_0$, 
and the Hilbert function of $F_1$
determines $F_1$, since any finitely generated 
free graded module is determined by its Hilbert function. 

Thus on \pr2 the problem of determining the modules in a resolution for $I_Z$ reduces to 
finding $h_{I_Z}$ and to computing for each $d$ the number $\nu_d$
of generators of degree $d$ in a minimal homogeneous set of generators for $I_Z$. 

Another feature of working on \pr2 is the ease with which one
can extend the notion of fat point subschemes supported at distinct points 
to include the possibility of infinitely near points. In fact, it is this extended notion
of fat point subschemes that we will use in this paper and which we now introduce.

We first put the notion of fat point subscheme into a context which will make
our extended notion natural. Let $p_1,\ldots,p_r\subset\pr 2$
be distinct points of the plane. Let $\pi:X\to\pr 2$ be the blowing up of 
each of the points. Then the divisor class group
$\hbox{Cl}(X)$ (elements of which, when convenient, we will identify 
with the corresponding associated invertible sheaves)  of $X$ is a free abelian
group of rank $r+1$, with basis $e_0,\ldots,e_r$, where $e_0$ is the pullback to $X$ 
of the class of a line, and $e_1,\ldots,e_r$ are the classes of the exceptional divisors
$E_1,\ldots,E_r$ of the blowings up of the points $p_1,\ldots,p_r$. Let $\C I_Z$
be the sheaf of ideals of some fat point subscheme $Z=m_1p_1+\cdots+m_rp_r$,
and let $\C F_n$ denote $ne_0-m_1e_1-\cdots-m_re_r$.
Then $\C I_Z=\pi_*(\C F_0)$, and $\hbox{dim}((I_Z)_d)=h^0(\pr2,\C I_Z\otimes\C O_{\pr2}(d))=
h^0(X,\C F_d)$ for every $d$ (so $I_Z$ is isomorphic as a graded module to 
$\bigoplus_{d} H^0(X,\C F_d)$), and in fact
$h^i(\pr2,\C I_Z\otimes\C O_{\pr2}(d))=h^i(X,\C F_d)$ for every $i$.

Moreover, $(I_Z)_d\otimes R_1\to(I_Z)_{d+1}$
corresponds to $H^0(X,\C F_d)\otimes H^0(X,e_0)\to H^0(X,\C F_{d+1})$,
and so the kernel and cokernel of the former have the same dimension as for the latter.
Given divisor classes \C G and \C H, we will denote the kernel and cokernel
of $H^0(X,\C G)\otimes H^0(X,\C H)\to H^0(X,\C G+\C H)$ by $\C R(\C G,\C H)$ and 
$\C S(\C G,\C H)$, respectively, and their dimensions by $\r(\C G,\C H)$ and $\s(\C G,\C H)$.
Thus $h_{I_Z}(d)=h^0(X,\C F_d)$ and $\nu_{d+1}(Z)=\s(\C F_d,e_0)$, so 
we see the Hilbert function of $I_Z$ and the number of generators of $I_Z$
in each degree can be found by studying invertible sheaves
on $X$. And indeed, this is the approach we take in this paper. 

To subsume the case of infinitely near points we now define the notion
of {\it essentially distinct} points. Let $p_1\in X_0=\pr2$,
and let $p_2\in X_1$, $\ldots$, $p_r\in X_{r-1}$, where, for $0\le i\le r-1$,
$\pi_i:X_{i+1}\to X_i$ is the blowing up of $p_{i+1}$. We will denote $X_r$ by $X$
and the composition $X\to \pr2$ by $\pi$. We call the indexed points
$p_1,\ldots,p_r$ {\it essentially distinct\/} points of \pr2; note that $p_j$
for $j>i$ may be infinitely near $p_i$. Denoting the class of the
1-dimensional scheme-theoretic fiber $E_i$
of $X_r\to X_i$ by $e_i$ and the pullback to $X_r$ of the class of a line in \pr2 by $e_0$,
we have, as in the case of distinct points, the basis $e_0,\ldots,e_r$
of the divisor class group of $X$ corresponding
to $p_1,\ldots,p_r$, and which we will call an {\it exceptional configuration}. 
Then $\pi_*(-m_1e_1-\ldots-m_re_r)$ is a coherent sheaf of ideals
on \pr2 defining a 0-dimensional subscheme $Z$ generalizing the usual notion of fat
point subscheme. In analogy with the notation used above, we will denote $Z$ by
$m_1p_1+\cdots+m_rp_r$ and refer to $Z$ as a fat point subscheme. 
As an aside we also note that the stalks of $\pi_*(-m_1e_1-\ldots-m_re_r)$ 
are complete ideals in the stalks of the local rings of the structure
sheaf of \pr2, and that conversely if \C I is a coherent sheaf of ideals on \pr2
whose stalks are complete ideals and if \C I defines a 0-dimensional subscheme,
then there are essentially distinct points $p_1,\ldots,p_r$ of \pr2
and integers $m_i$ such that with respect to the corresponding
exceptional configuration we have $\C I=\pi_*(-m_1e_1-\cdots-m_re_r)$.
Thus our extended notion of fat points is precisely what is obtained
by considering 0-dimensional subschemes defined by coherent sheaves of ideals
whose stalks are complete ideals.

Allowing the possibility of infinitely near points necessitates
dealing with certain technicalities.
In particular, the subscheme $Z$ does not uniquely determine 
$-m_1e_1-\cdots-m_re_r$. For example, if $p_1$ and $p_2$ are distinct
points of \pr2, then $\pi_*(-e_1+e_2)=\pi_*(-e_1)$ both give the sheaf of ideals
defining the subscheme $Z=p_1$. To get uniqueness, we recall that 
the divisor class group of $X$ supports an intersection form, with respect to which
the exceptional configuration $e_0,\ldots,e_r$ is orthogonal with
$-1=-e_0^2=e_1^2=\cdots=e_r^2$. The inequalities $(-m_1e_1-\cdots-m_re_r)\cdot C_{ij}\ge 0$,
where the index $i$ runs over the divisors $E_i$ and 
$j$ runs over the components $C_{ij}$ of $E_i$, corresponds to what older
terminology called the {\it proximity inequalities}. Thus we will say that
a divisor class $\C F$ on $X$ satisfies the proximity inequalities if
$\C F\cdot C\ge 0$ for every component $C$ of each divisor $E_i$. Moreover,
given essentially distinct points $p_1,\ldots,p_r$ and a subscheme
$Z=m_1p_1+\cdots+m_rp_r$, we will abbreviate saying that the class 
$-m_1e_1-\cdots-m_re_r$ coming from the coefficients 
$m_1,\ldots,m_r$ used to define $Z$ satisfies the proximity inequalities
by simply saying that $m_1p_1+\cdots+m_rp_r$ satisfies the proximity inequalities.
Uniqueness is now a consequence of the fact that 
if $\pi_*(-a_1e_1-\cdots-a_re_r)=\pi_*(-b_1e_1-\cdots-b_re_r)$,
where $-a_1e_1-\cdots-a_re_r$ and $-b_1e_1-\cdots-b_re_r$ both satisfy
the proximity inequalities, then $a_i=b_i$ for each $i$. In particular,
we have a bijection between subschemes of fat points in \pr2
and 0-cycles $m_1p_1+\cdots+m_rp_r$, where $p_1,\ldots,p_r$ are essentially
distinct points of \pr2 and $m_1p_1+\cdots+m_rp_r$ satisfies the
proximity inequalities. 

From another point of view, the significance of the proximity inequalities
is given by an old and well known result 
saying that the linear system of sections of $de_0-m_1e_1-\cdots-m_re_r$ is fixed component free
for $d$ sufficiently large if (and only if) $-m_1e_1-\cdots-m_re_r$
satisfies the proximity inequalities. 
The proximity inequalities also manifest themselves even in
the usual case that $p_1,\ldots,p_r$ are distinct points of \pr2, 
since in this case $m_1p_1+\cdots+m_rp_r$ satisfying the proximity
inequalities just means that each coefficient $m_i$ is nonnegative,
which is generally taken for granted without comment.

We now discuss in more detail the approach we take in this paper.
Let $Z=m_1p_1+\cdots+m_rp_r$ be a fat point subscheme of \pr2, satisfying the proximity
inequalities. Let $e_0,\ldots,e_r$ be the exceptional configuration corresponding
to the essentially distinct points $p_1,\ldots,p_r$. For each $d$, let $\C F_d$
be the class $de_0-m_1e_1-\cdots-m_re_r$. As we have seen above, determining the graded modules
in a minimal free resolution of $I_Z$ amounts to computing $h^0(X,\C F_d)$ and
$\s(\C F_d,e_0)$ for each $d$. If $\C F_d$ is not the class of an effective divisor,
then clearly $h^0(X,\C F_d)=0$ and $\s(\C F_d,e_0)=h^0(X,\C F_{d+1})$.
If, however, $\C F_d$ is the class of an effective divisor, then $\C F_d$ has
a {\it Zariski decomposition\/}: $\C F_d=\C H+\C N$, where $h^0(X,\C N)=1$
and \C H is {\it numerically effective} (i.e., meets every effective divisor
nonnegatively) with $h^0(X,\C F_d)=h^0(X,\C H)$. By \ir{cokfact}, 
$\s(\C F_d,e_0)=\s(\C H,e_0)+h^0(X,\C F_{d+1})-h^0(X,\C H+e_0)$. Thus,
to implicitly determine a resolution of $I_Z$, it is enough to determine: the
monoid EFF of divisor classes of effective divisors; a Zariski decomposition for each
class in EFF; and $h^0(X, \C H)$ and $\s(\C H,e_0)$
for each class \C H in the cone NEFF of numerically effective classes.

If $X$ is any smooth projective rational surface for which $-K_X$ is effective,
then EFF and NEFF can be found, as can a Zariski decomposition for any class in EFF
and $h^0(X, \C H)$ for any numerically effective class \C H, by applying the results
of [\ars]. In particular, since in terms of an exceptional configuration $e_0,\dots,e_r$
on a blowing up $X$ of \pr2 at essentially distinct points $p_1,\ldots,p_r$ we always have
$-K_X=3e_0-e_1-\cdots-e_r$, we see $-K_X$ is the class of an effective divisor whenever 
$p_1,\ldots,p_r$ lie on a curve of degree 3 or less.
Thus the novel part in determining resolutions for fat point subschemes of \pr2
supported at points of a cubic is in determining $\s(\C H,e_0)$
for numerically effective classes on the blowing up $X$ of the points.
What we will see here is that if $X$ is the blowing up of
points on a conic, $\s(\C H,e_0)=0$ for
any numerically effective class \C H. This is no longer true
for points on a cubic, which is partly why our results for the cubic case
are less comprehensive.

\irrnSection{Background on Surfaces}{bkgdsect}
This section recalls results on surfaces that we will need
later. For those results which are not standard or well known
we give an indication of proof; for the reader's convenience at the least we 
usually provide
references for the rest. Given a subvariety $C\in X$ and a class \C L
on $X$, it will be convenient, if our meaning is clear, to write
$H^i(C, \C L)$ for the cohomology of the restriction, rather than 
$H^i(C, \C L\otimes\C O_C)$.

\prclm{Lemma}{one}{Let $\pi:Y\to X$ be a birational morphism of 
smooth projective rational surfaces,
$\pi^*:\hbox{Cl}(X)\to\hbox{Cl}(Y)$ the corresponding homomorphism 
on divisor class groups, and let \C L be a divisor class on $X$.
Then $\pi^*$ is an injective intersection-form preserving
map of free abelian groups of finite rank; there is a natural isomorphism
$H^i(X,\C L)=H^i(Y,\pi^*\C L)$ for every $i$; and \C L is the class of an effective
divisor (resp., numerically effective) if and only if $\pi^*\C L$ is.}

\Prf See [\ars] for an indication of proof. \qed

By \ir{one}, little harm is done by identifying $\hbox{Cl}(X)$ with its image
in $\hbox{Cl}(Y)$. This is also compatible with exceptional configurations.
For suppose that $X=X_r$ is obtained by blowing up essentially distinct points
$p_1,\ldots,p_r$ of $\pr2=X_0$, where, for $1\le i\le r$, $X_i$ is the blowing up of
$p_1,\ldots,p_i$. If $\pi:X_j\to X_i$ is, for some $i<j\le r$, the blowing up of
$p_{i+1},\ldots,p_j$, with $e_0,\ldots,e_i$ the exceptional configuration on $X_i$
corresponding to $p_1,\ldots,p_i$ and $e_0',\ldots,e_j'$ the exceptional configuration on $X_j$
corresponding to $p_1,\ldots,p_j$, then $\pi^*(e_l)=e_l'$ for $l\le i$. For simplicity
then, we will for each $0\le l\le r$ simply denote the exceptional configuration
on $X_l$ corresponding to $p_1,\ldots,p_l$ by $e_0,\ldots,e_l$, and leave to context
which surface $X_j$, $j\ge i$ we wish at any given time to regard $e_i$ as a class on.

We now recall some facts, the most important of which
is (a), the formula of Riemann-Roch for a rational surface ([14]).
For proofs of (b) and (c), we refer to [\mtnwest].

\prclm{Lemma}{RR}{Let $X$ be a smooth projective 
rational surface, and let \C F be
a divisor class on $X$.
\item{(a)} We have: $h^0(X, \C F)-h^1(X, \C F) + h^2(X, 
\C F)= (\C F^2-K_X\cdot \C F)/2+1$.
\item{(b)} If \C F is the class of an effective divisor, then $h^2(X, \C F)=0$.
\item{(c)} If \C F is numerically effective,
then $h^2(X, \C F)=0$ and $\C F^2\ge 0$.}

Here we recall some results from [18], where, in line with the notation in
\ir{introsect} and following [18], given coherent sheaves
\C A and \C B on a scheme $T$, we will denote the cokernel of
the natural map $H^0(X, \C A)\otimes H^0(X, \C B)\to H^0(X, \C A\otimes\C B)$
by $\C S(\C A,\C B)$ and the kernel by $\C R(\C A,\C B)$. Also,
$\Gamma$ denotes the global sections functor.

\prclm{Proposition}{Mumford}{Let $T$ be a closed subscheme of projective space,
let \C A and \C B be coherent sheaves on $T$ and let \C C be the class of an
effective divisor $C$ on $T$.
\item{(a)} If the restriction homomorphisms $H^0(T, \C A)\to H^0(C,\C A\otimes\C O_C)$
and $H^0(T, \C A\otimes\C B)\to H^0(C,\C A\otimes\C B\otimes\C O_C)$ are surjective
(for example, if $h^1(T,\C A\otimes\C C^{-1})=0=h^1(T,\C A\otimes\C C^{-1}\otimes\C B)$),
then mapping the terms of the exact sequence 
$(0\to \Gamma(\C A\otimes\C C^{-1}) \to
\Gamma(\C A)\to\Gamma(\C A\otimes\C O_C)\to 0)\otimes\Gamma(\C B)$
to those of $0\to \Gamma(\C A\otimes\C C^{-1}\otimes \C B) \to
\Gamma(\C A\otimes \C B)\to\Gamma((\C A\otimes \C B)\otimes\C O_C)\to 0$
leads to the exact sequence
$$\eqalign{0\to & \C R(\C A\otimes\C C^{-1},\C B)\to
\C R(\C A,\C B)\to\C R(\C A\otimes\C O_C,\C B)\to\cr
&\C S(\C A\otimes\C C^{-1},\C B)\to\C S(\C A,\C B)\to\C S(\C A\otimes\C O_C,\C B)\to0.\cr}$$
\item{(b)} If $H^0(T, \C B)\to H^0(C,\C B\otimes\C O_C)$ is surjective
(for example, if $h^1(T,\C B\otimes\C C^{-1})=0$), then 
$\C S(\C A\otimes\C O_C,\C B)=\C S(\C A\otimes\C O_C,\C B\otimes\C O_C)$.
\item{(c)} If $T$ is a smooth curve of genus $g$, and \C A and \C B are line 
bundles of degrees at least $2g+1$ and $2g$, respectively, then 
$\C S(\C A,\C B)=0$.}

\Prf See [18] for (a) and (c); we leave (b) as an easy exercise for the reader.\qed

It will be helpful to generalize \ir{Mumford}(c) to nonsmooth curves with $g=0$.
We do so in \ir{keylemma}, using the following technical result, proved in [\ars].

\prclm{Lemma}{keycor}{Let $X$ be a smooth projective 
rational surface
and let \C N be the class of a nontrivial effective divisor $N$ 
on $X$. If $\C N+K_X$ is
not the class of an effective divisor and \C F meets every component of \C N nonnegatively, 
then $h^0(N, \C F)>0$, $h^1(N, \C F)=0$, $N^2+N\cdot K_X<-1$,
and every component $M$ of $N$ is a smooth rational curve (of negative
self-intersection, if $M$ does not move).}

\prclm{Lemma}{keylemma}{Let $X$ be a smooth projective 
rational surface, and let \C N be the class of an effective divisor $N$ on $X$
such that $h^0(X, \C N+K_X)=0$. If \C F and \C G are the restrictions to $N$
of divisor classes $\C F'$ and $\C G'$ on $X$ which meet each component of $N$ nonnegatively,
then $\C S(\C F,\C G)=0$.}

\Prf To prove the lemma, induct on the sum $n$ of 
the multiplicities of the components of $N$.
By \ir{keycor}, $h^1(N, \C O)=0$ and every component
of $N$ is a smooth rational curve. Thus the case $n=1$
is trivial (since then $N=\pr 1$, and
the space of polynomials of degree $f$ in two variables 
tensor the space of polynomials of degree $g$ in two variables 
maps onto the space of polynomials of degree $f+g$). So say $n>1$.

As in the proof
of Theorem 1.7 of [1] (or see the proof of
Lemma II.6 of [\mtnwest]), $N$ has a component $C$ such that 
$(N-C)\cdot C\le 1$. Let $L$ be the effective divisor $N-C$
and let \C L be its class.
Thus we have an exact sequence 
$0\to \C O_C\otimes(-\C L)\to \C O_N\to \C O_L\to 0$.
Now, $-L\cdot C\ge -1$, and both 
$\C F'$ and $\C G'$ meet $C$ nonnegatively. We may assume
$\C F'\cdot C\ge\C G'\cdot C$, otherwise reverse the roles of 
$\C F'$ and $\C G'$. Since
$C=\pr 1$, we see that $h^1(C, \C O_C\otimes(\C F'-\C L))$,
$h^1(C, \C O_C\otimes(\C G'-\C L))$ and 
$h^1(C, \C O_C\otimes(\C F'+\C G'-\C L))$ all vanish.
An argument similar to that used to prove \ir{Mumford}(a, b) 
now shows that we have an exact sequence
$\C S(\C O_C\otimes(\C F'-\C L),\C O_C\otimes\C G')\to
\C S(\C F,\C G)\to\C S(\C O_L\otimes\C F,\C O_L\otimes\C G)\to0$.
Since $\C S(\C O_L\otimes\C F,\C O_L\otimes\C G)=0$ by induction,
it suffices to show
$\C S(\C O_C\otimes(\C F'-\C L),\C O_C\otimes\C G')=0$.
If $C\cdot(\C F'-\C L)\ge 0$, then the latter is 0 (as in the previous paragraph).
Otherwise, we must have $0=\C F'\cdot C=\C G'\cdot C$ and
$C\cdot L=1$, so
$\C O_C(-1)=\C O_C\otimes(\C F'-\C L)$ and
$\C O_C=\C O_C\otimes\C G'$, which means 
$h^0(\C O_C,\C O_C\otimes(\C F'+\C G'-\C L))=0$ 
and hence again $\C S(\C O_C\otimes(\C F'-\C L),\C O_C\otimes\C G')=0$. \qed

When working with distinct points $p_1,\ldots,p_r$ it can be convenient
to reindex them; for example, so that an expression
$m_1p_1+\cdots+m_rp_r$ has $m_1\ge \cdots\ge m_r$.
If $p_1,\ldots,p_r$ are only essentially distinct, reindexing in this
way, properly speaking, makes no sense (since $p_i$ lives on the surface
obtained by blowing up $p_1,\ldots,p_{i-1}$). Nonetheless, a reindexation 
that preserves the natural partial ordering of ``infinite nearness''
seems intuitively acceptable. We make a short digression to justify
this intuition.

Suppose $p_1,\ldots,p_r$ and $q_1,\ldots,q_r$ are essentially
distinct points of \pr2. Let $X$ be the blowing up of 
$p_1,\ldots,p_r$ with $e_0,\ldots,e_r$ being the associated exceptional configuration,
and let $X'$ be the blowing up of 
$q_1,\ldots,q_r$ with $e_0',\ldots,e_r'$ being the associated exceptional configuration.
If there is an isomorphism $f:X\to X'$ such that $f^*(e_0')=e_0$, then there is
a unique permutation $\sigma_f$  of $\{1,\ldots,r\}$ such that 
$f^*(e_i')=e_{\sigma_f(i)}$ for every $i\ge 1$, and
it follows that any subscheme $Z=m_1p_1+\cdots+m_rp_r$ is projectively
equivalent to $Z'=m_{\sigma_f(1)}q_1+\cdots+m_{\sigma_f(r)}q_r$. 
Thus we shall say that a bijection $\sigma:\{q_1,\ldots,q_r\}\to\{p_1,\ldots,p_r\}$ is
an {\it equivalence\/} and that $p_1,\ldots,p_r$ and $q_1,\ldots,q_r$ are {\it equivalent \/}
if for some $f$ as above we have $\sigma(q_i)=p_{\sigma_f(i)}$ for every $i$.
Similarly, we shall say that a permutation of $\{1,\ldots,r\}$
is an {\it equivalence\/} if it is $\sigma_f$ for some such $f$.

By the next lemma
we see that we may always assume that a subscheme
$Z=m_1p_1+\cdots+m_rp_r$, satisfying the proximity inequalities,
also satisfies $m_1\ge \cdots\ge m_r$ up to equivalence.

\prclm{Lemma}{reindex}{Let $p_1,\ldots,p_r$ be essentially distinct points of \pr2.
\item{(a)} Any permutation $\sigma$ of $\{1,\ldots,r\}$, such that 
$\sigma(j)\ge \sigma(i)$ whenever $p_j$ is infinitely near $p_i$,
is an equivalence.
\item{(b)} If $Z=m_1p_1+\cdots+m_rp_r$ satisfies the proximity inequalities,
then $m_j\le m_i$ whenever $p_j$ is infinitely near to $p_i$;
in particular, up to equivalence we have $m_1\ge \cdots\ge m_r$.}

\Prf Let $X$ be the blowing up
of $p_1,\ldots,p_r$ , with $e_0,\ldots,e_r$  being the corresponding exceptional
configuration.

(a) One merely needs to check that $e_0,e_1',\ldots,e_r'$ is also an exceptional
configuration, where $e_i'=e_{\sigma(i)}$ for each $i>0$. But this
follows from Theorem 1.1 of [\duke].

(b) Since $p_j$ being infinitely near to $p_i$ implies $e_i-e_j$ is the 
class of an effective divisor,
the proximity inequalities imply $-(m_1e_1+\cdots+m_re_r)\cdot (e_i-e_j)\ge 0$
and hence $m_i\ge m_j$. In particular, we can choose $l$ such 
that $m_l\le m_i$ for all $1\le i\le r$ and such that
only $p_l$ among $p_1,\ldots,p_r$ is infinitely near to $p_l$; then the permutation 
$\sigma$ which is the identity on $\{1,\ldots,l-1\}$, and for which
$\sigma(l)=r$ and $\sigma(i)=i-1$ for $l<i\le r$, is an equivalence by (a).
I.e., up to equivalence we may assume $m_r$ is least among $m_1,\ldots,m_r$,
and the result follows by induction.   \qed

We now have a lemma well known for distinct points. As a consequence
of \ir{keycor}, it holds more generally for essentially distinct points too.

\prclm{Lemma}{sg}{Let $p_1,\ldots,p_r$ be essentially distinct points
of \pr2 and let $e_0, e_1,\ldots,e_r$ be the corresponding
exceptional configuration on the blowing up $X$ of \pr2 at the points. 
Suppose that $m_1,\ldots,m_r$ are nonnegative integers
such that $-(m_1e_1+\cdots+m_re_r)$ satisfies the proximity inequalities. 
Then $h^0(X, \C F_d)>0$ and $h^1(X, \C F_d)=0$ for $d\ge -1+(m_1+\cdots+m_r)$, where
$\C F_d=de_0-(m_1e_1+\cdots+m_re_r)$.}

\Prf Suppose there is at most one index $i$ with $m_i>0$.
By \ir{reindex} we may assume $r=1$, in which case we just
need to check that $h^1(X, de_0-m_1e_1)=0$
for $d\ge -1+m_1$, which is straightforward, by restricting 
to a general section $B$ of $e_0-e_1$.

We now consider $\C F_d$ in case $m_i>0$ for at least two indices $i$.
By \ir{one} we may assume that $m_r>0$.
Then $-(m_1e_1+\cdots+(m_r-1)e_r)$ satisfies the proximity inequalities,
so by induction we may assume that $h^1(X, (d-1)e_0-(m_1e_1+\cdots+(m_r-1)e_r))=0$.
Also note that $(d-1)e_0-(m_1e_1+\cdots+(m_r-1)e_r)=\C F_d-(e_0-e_r)$.
Clearly $e_0-e_r$ is the class of an effective divisor; $p_r$ is infinitely 
near a bona fide point of \pr2,
and that point is $p_i$ for some $i\le r$. Then $e_0-e_r=(e_0-e_i)+
(e_i-e_r)$; $e_0-e_i$ is numerically effective, its linear system 
of sections being the pencil of lines through the point $p_i$, 
while $E_r$ is a component of $E_i$, and $e_i-e_r$ is the class
of the difference $E_i-E_r$. Thus we may choose a section $B$ of $e_0-e_r$
whose components, apart from a numerically effective divisor,
are all components of $E_i$. If $C$ is a component of $E_i$, then 
$F_d\cdot C=-(m_1e_1+\cdots+m_re_r)\cdot C\ge 0$ by hypothesis.
Since any numerically effective divisor meets any effective
divisor nonegatively, it will follow that $\C F_d$ meets every component
of $B$ nonnegatively if we show $\C F_d$ is the class of an effective divisor. 
But $d+1\ge m_1+\cdots+m_r$, the right hand side 
of which involves at least
two terms; thus $(d+1)^2>m_1^2+\cdots+m_r^2$, hence $\C F_d^2-K_X\cdot \C F_d=
-1+[(d+1)^2-(m_1^2+\cdots+m_r^2)] +[d-(m_1+\cdots+m_r)]\ge -1$,
so by \ir{RR}, $h^0(X, \C F_d)>0$ (since $d\ge 0$ implies 
$h^2(X, \C F_d)=0$ by duality). 

Next, $e_0\cdot(K_X+(e_0-e_r))<0$ (because $e_0\cdot K_X=-3$), 
so $h^0(X, K_X+(e_0-e_r))=0$, since $e_0$ is numerically effective. 
We now can apply \ir{keycor} to obtain $h^1(B, \C F_d\otimes\C O_B)=0$, from
which $h^1(X, \C F_d)=0$ follows by taking cohomology
of $0\to \C F_d-(e_0-e_r)\to \C F_d\to \C F_d\otimes\C O_B\to 0$.
\qed

The next result concerns vanishing of $\C S (\C F,\C G)$.

\prclm{Theorem}{cokvan}{Let $X$ be the blowing up of essentially distinct points
$p_1,\ldots,p_r$ of \pr2, and let $e_0,\ldots,e_r$ be the corresponding 
exceptional configuration. Suppose $\C F=a_0e_0-a_1e_1-\cdots-a_re_r$
and $\C G=b_0e_0-b_1e_1-\cdots-b_re_r$ satisfy the proximity inequalities
and that $a_0\ge a_1+\cdots+a_r$ and $b_0\ge b_1+\cdots+b_r$.
Then $\C S (\C F,\C G)=0$.}

\Prf The result is true and easy to see if \C F and \C G are multiples of $e_0$,
so assume that $a_r$ say is positive. Then $\C H=\C F-(e_0-e_r)$ satisfies the
proximity inequalities and also has $e_0\cdot\C H\ge \sum_{i>0}e_i\cdot\C H$,
so by induction we may assume that $\C S(\C H,\C G)=0$. By \ir{sg}, we have
$h^1(X,\C H)=0=h^1(X,\C H+\C G)$. If $b_r>0$, then \ir{sg} also gives $h^1(X,\C G-(e_0-e_r))=0$. If 
$b_r=0$, then $h^1(X,\C G-e_0)=0$ by \ir{sg} and $e_r$ is a fixed component of
$\C G-(e_0-e_r)$, so $h^0(X, \C G-(e_0-e_r))=h^0(X, \C G-e_0)$. This, together with 
$(\C G-(e_0-e_r))^2-K_X\cdot(\C G-(e_0-e_r))=(\C G-e_0)^2-K_X\cdot(\C G-e_0)$ and
\ir{RR} implies $h^1(X,\C G-(e_0-e_r))=h^1(X,\C G-e_0)=0$. Now from \ir{Mumford}(a,b)
we have an exact sequence 
$0=\C S(\C H,\C G)\to \C S(\C F,\C G)\to \C S(\C F\otimes\C O_L,\C G\otimes\C O_L)\to 0$,
where $L$ is a general section of $e_0-e_r$. 
As in the proof of \ir{sg}, $K_X+(e_0-e_r)$ is not the class of an effective divisor
and \C F and \C G meet each component of $L$
nonnegatively, so $\C S(\C F\otimes\C O_L,\C G\otimes\C O_L)=0$ by \ir{keylemma},
and hence $\C S(\C F,\C G)=0$, as required.   \qed

The following result, giving another vanishing 
criterion, is well known (see Proposition 3.7 of [4]) and follows easily
by appropriately applying \ir{Mumford}; it essentially says that
no generator of $I_Z$ need be taken in degrees greater than the regularity
of $I_Z$.

\prclm{Lemma}{myomega}{Let $e_0,\ldots,e_r$ be the exceptional configuration
corresponding to a blowing up $X\to\pr2$ at essentially
distinct points $p_1,\ldots, p_r$. Let $Z=m_1p_1+\cdots+m_rp_r$
satisfy the proximity inequalities,
and let $\C F_d$ denote $de_0-m_1e_1-\cdots-m_re_r$.
If $h^1(X,\C F_t)=0$, then $\C S(\C F_d,e_0)=0$ for $d>t$.}

In the following lemma, given the class \C F of an effective divisor $F$
with $N$ denoting the fixed components of the complete linear system $|F|$,
we denote the class of $N$ 
by $(\C F)_f$ (this is the {\it fixed\/} 
part of \C F) and $\C F-(\C F)_f$ by $(\C F)_m$ (the free or {\it moving\/} part). 
Thus $\C F=(\C F)_m+(\C F)_f$
is a Zariski decomposition of \C F.

\prclm{Lemma}{cokfact}{Let $e_0,\ldots,e_r$ be an exceptional configuration
on a surface $X$ and let \C F be a class on $X$.
\item{(a)} If \C F is not the class of an effective divisor,
then $\s(\C F,e_0)=h^0(X,\C F+e_0)$.
\item{(b)} If \C F is the class of an effective divisor with a Zariski decomposition 
$\C F=\C H+\C N$ (with \C H being the numerically effective part),
then $\s(\C F,e_0)=\s(\C H,e_0) + h^0(X,\C F+e_0)-h^0(X,\C H+e_0)$.
\item{(c)} If \C F is the class of an effective divisor and $(\C F)_f\ne (\C F+e_0)_f$,
then $\s(\C F,e_0)\ge h^0(X,\C F+e_0)-h^0(X,(\C F)_m+e_0)>0$.}

\Prf (a) This is clear so consider (b). Regarding $\C H$ and \C F as sheaves, we have an inclusion
$\C H\to\C F$ which induces an isomorphism
on global sections. Thus we have a commutative diagram with exact columns

$$\hbox{\hfil\vbox{\halign{\hfil $#$\hfil &\hfil $#$\hfil &\hfil $#$\hfil \cr
       0                          &     &           0            \cr
     \downarrow                   &     &      \downarrow        \cr
H^0(X,\C H)\otimes H^0(X,e_0) & \to & H^0(X,\C H+e_0)    \cr
     \downarrow                   &     &      \downarrow        \cr
H^0(X,\C F)\otimes H^0(X,e_0)     & \to & H^0(X,\C F+e_0)        \cr
     \downarrow                   &     &                        \cr
       0                          &     &                        \cr}}\hfil}$$

\noindent and the image of $H^0(X,\C H)\otimes H^0(X,e_0) \to H^0(X,\C F+e_0)$
equals the image of 
$H^0(X,\C F)\otimes H^0(X,e_0)\to H^0(X,\C F+e_0)$, which means
that $\s(\C F,e_0)=\s(\C H,e_0) + h^0(X,\C F+e_0)-h^0(X,\C H+e_0)$.

(c) By (b) we know $\s(\C F,e_0)\ge h^0(X,\C F+e_0)-h^0(X,(\C F)_m+e_0)$.
But $h^0(X,\C F+e_0)-h^0(X,(\C F)_m+e_0)$ just measures the extent to which
there are sections of $\C F+e_0$ not containing the full fixed part of \C F.
Since  $(\C F)_f\ne (\C F+e_0)_f$, this is positive. \qed

\irrnSection{Resolutions}{resolsect}
In this section we will, under certain restrictions, 
study minimal free resolutions of homogeneous ideals $I_Z$
(over the homogeneous coordinate ring $k[\pr2]$ of \pr2,
which hereafter we will denote by $R$),
where $p_1,\ldots,p_r$ are essentially distinct points of \pr2
and $Z=m_1p_1+\cdots+m_rp_r$ satisfies the proximity inequalities.
By \ir{reindex}, we may, if it is convenient, assume that 
$m_1\ge\cdots\ge m_r>0$. 

Given essentially distinct points $p_1,\ldots,p_r$ of \pr2, we say that
the points are points of a curve of degree $n$ if, on the blowing up $X$ of \pr2
at the points, $ne_0-e_1-\cdots -e_r$ is the class of an effective divisor. We say
that $p_1,\ldots,p_r$ are points of a curve of degree $n$ with some property 
(say smooth or irreducible, etc.) if $ne_0-e_1-\cdots -e_r$ is the class
of such a curve. Our results involve essentially distinct points $p_1,\ldots,p_r$ on
curves in \pr2 of degree at most 3. Our results allow one
recursively to determine the graded modules in a minimal free resolution
of $I_Z$ for any fat point subscheme $Z=m_1p_1+\cdots+m_rp_r$ 
satisfying the proximity inequalities, as long as $p_1,\ldots,p_r$
are points of a conic, or, in certain cases, points of a smooth cubic.
Since essentially distinct points on a line certainly are points
on a conic, our results apply also to points on a line; 
the case of points on a line
turns out to be simple enough to write down the resolution 
completely explicitly, which we do in \ir{linecase}.

\irSubsection{Points on a conic}{conicsubsect}

We now consider points, possibly infinitely near, on a conic, possibly nonsmooth.
If $X$ is the blowing up of \pr2 at such points, we show $\s(\C F,e_0)=0$
for any $\C F\in\hbox{NEFF}$. As demonstrated in \ir{conicexample}, this allows us to work out
resolutions in any specific case of points on a conic. Also, in the special case of points on a line
(which is subsumed by the case of points on a conic), we give a completely general and explicit result
in \ir{linecase}.

We begin with a lemma: part (a) will be used in the proof of \ir{cokonaconic};
the other parts will be helpful references in working out examples, such as \ir{conicexample}.

\prclm{Lemma}{eff}{Let $X$ be the blowing up of essentially distinct points $p_1,\ldots,p_r$ 
of a conic in \pr2; i.e., $2e_0-e_1-\cdots-e_r$ is the class of an effective
divisor, where $e_0,\ldots,e_r$ is the exceptional configuration
corresponding to $p_1,\ldots,p_r$. Also, for each $i>0$, recall 
$E_i$ denotes the unique effective divisor whose class is $e_i$.

\item{(a)} If \C C is the class of a reduced and irreducible curve 
of negative self-intersection, then \C C is either:
the class of a component of $E_i$, $i>0$; the class
of a component of an effective divisor $Q$ with $\C C\cdot Q<0$, where
$\C Q=2e_0-e_1-\cdots-e_r$ is the class of $Q$; or the class
$e_0-e_i-e_j$, for some $0<i<j$.

\item{(b)} If $\C F\in\hbox{NEFF}$, then $\C F\in\hbox{EFF}$, and
\C F is regular (i.e., $h^1(X,\C F)=0$) and its linear system of sections is base point
(and hence fixed component) free.

\item{(c)} If $r\ge 2$, then:
\itemitem{(i)} NEFF consists of the classes $\C F\in
\hbox{Cl}(X)$ such that $\C F\cdot \C C\ge 0$ whenever
\C C is the class of an irreducible curve of negative self-intersection, and
\itemitem{(ii)} EFF is generated by the classes of curves of 
negative self-intersection.}

\Prf (a) Let \C C be the class of a reduced and irreducible curve $C$ of arithmetic genus $g$
and of negative self-intersection. Since $C$ is effective and
$e_0$ is numerically effective, we have $C\cdot e_0\ge 0$. 

If $C\cdot e_0=0$, then $C$ is a component of some $E_i$, $i>0$.

If $C\cdot e_0=1$, then \C C must be of the form
$e_0-e_{i_1}-\cdots-e_{i_l}$, with $0<i_1<\cdots<i_l$. Since
$C^2<0$, we have $l\ge 2$; if $l>2$ then $C$ meets \C Q negatively
and hence is a component of $Q$. 

Say $C\cdot e_0=2$. Then $C\cdot e_i\ge 0$ for all $i>0$, since $C$ is clearly not a component
of any $E_i$. Also, since any reduced and irreducible plane conic
is smooth, we see $C\cdot e_i\le 1$, $i>0$. Thus \C C is the sum
of \C Q with those $e_i$ with $C\cdot e_i=0$, so we see 
that \C C has a section coming from \C Q and those $e_i$ with $C\cdot e_i=0$,
and hence either that $C$ is not fixed
(contradicting $C^2<0$), or that $C\cdot e_i=1$ for all $i$
and therefore that $\C C=\C Q$ and $C\cdot\C Q=C^2<0$.

Suppose that $C\cdot e_0>2$. Since 
$-K_X=e_0+\C Q$ and clearly $C$ is not a component of \C Q, we
have $-C\cdot K_X\ge C\cdot e_0\ge 3$.
Hence by the adjunction formula we have
$0>C^2=2g-2-C\cdot K_X\ge 2g+1>0$, contradiction.

(b) Since $-K_X$ is the class of an effective divisor, $\hbox{NEFF}\subset\hbox{EFF}$
[\mtnwest]. For the rest, the case $\C F=0$ is clear, so assume \C F is
not trivial. By \ir{RR}(c), 
we have $\C F^2\ge 0$. Since the space $e_0^\perp\subset
\hbox{Cl}(X)$ of classes perpendicular to $e_0$
is negative definite, we see $\C F\cdot e_0>0$.
Together, this means either that $\C F\cdot e_0>1$ (in which case 
$\C F\cdot (-K_X)=\C F\cdot (\C Q+e_0)>1$ and our result
follows by [\ars]) or that $\C F\cdot e_0=1$ and hence 
$\C F$ is either of the form $e_0$ or $e_0-e_i$, $i>0$
(and again we have $\C F\cdot (-K_X)>1$ and our result follows by
[\ars]).

(c) Clearly, (i) follows from (ii). To prove (ii), 
let \C F be in EFF. Then subtracting off fixed components,
which by (b) must be curves of negative self-intersection,
we obtain the free part $(\C F)_m$ of \C F, which is
numerically effective. Thus $(\C F)_m$ is in the cone of classes which meet
\C Q, each $e_i$, and each $e_0-e_i-e_j$, $0<i<j$ nonnegatively.
By Proposition 1.5.3 [\vanc], this cone 
is generated by classes of the form $e_0$, $e_0-e_{i_1}$, 
$2e_0-e_{i_1}-e_{i_2}-e_{i_3}$, $3e_0-e_{i_1}-e_{i_2}-e_{i_3}-e_{i_4}-e_{i_5}-e_{i_6}$,
and $de_0-(d-1)e_{i_1}-e_{i_2}-\cdots -e_{i_{d+2}}$, $d\ge 2$, where in each expression the indices
$i_1,i_2,\cdots$ are nonzero and distinct. Now it is enough to check that each of these classes is a sum of
classes of curves of negative self-intersection, which is straightforward. \qed

\prclm{Theorem}{cokonaconic}{Let $X$ be the blowing up of essentially 
distinct points $p_1,\ldots,p_r$ of \pr2 such that
$\C Q=2e_0-e_1-\cdots-e_r$ is the class of an effective divisor, where $e_0,\ldots,e_r$
is the exceptional configuration corresponding to $p_1,\ldots,p_r$. Then
$\C S(\C F,e_0)=0$ for any numerically effective class \C F.}

\Prf We will induct on $e_0\cdot \C F$. If $e_0\cdot \C F=0$ then
$\C F=0$, and $\C S(\C F,e_0)=0$ is clear. Also, the case that $r<2$ is 
covered by \ir{cokvan}, so we may assume $r\ge 2$, and by reindexing
(see \ir{reindex}) we may assume $\C F\cdot e_1\ge \cdots \ge \C F\cdot e_r\ge 0$.
Since, if $\C F\cdot e_r=0$ we may as well just work on the blowing up
of \pr2 at $p_1,\dots,p_{r-1}$, we may assume in fact that $\C F\cdot e_i>0$ for all $i>0$.
Now by explicitly checking against the cases enumerated in \ir{eff}(a),
given the class \C G of any reduced and irreducible effective 
divisor of negative self-intersection, we see $(\C F-\C Q)\cdot \C G\ge 0$. Thus
$\C F-\C Q$ is numerically effective.

Since $K_X=-3e_0+e_1+\cdots+e_r$, we have 
$\C Q=-K_X-e_0$. Clearly, given any numerically effective class \C H,
$(\C H+e_0)^2>0$, so, by duality and Ramanujan vanishing 
(see the first paragraph of [19, Theorem, p.\ 121], which holds
in all characteristics),
$h^1(X,\C H-\C Q)=h^1(X,-(\C H+e_0))=0$. Thus, taking $Q$ to be an 
effective divisor in the class \C Q, the exact sequence
$0\to \C H-\C Q\to \C H\to \C H\otimes\C O_Q\to 0$
is exact on global sections. In particular, this follows 
taking \C H to be either $e_0$, \C F or $\C F+e_0$.
In the former case we have
$\C S(\C F\otimes\C O_Q,e_0)=\C S(\C F\otimes\C O_Q,e_0\otimes\C O_Q)$
by \ir{Mumford}(b), and the latter two cases show that \ir{Mumford}(a) applies.
Since $\C S(\C F\otimes\C O_Q, e_0\otimes\C O_Q)$ vanishes by
\ir{keylemma}, and we may, by induction, assume $\C S(\C F-\C Q,e_0)=0$,
our result, $\C S(\C F,e_0)=0$, now follows from the exact sequence of
\ir{Mumford}(a).  \qed

\rem{Example}{conicexample}
We now give an example showing how to apply the results above to
work out resolutions. Let $L_1$ and $L_2$ be distinct lines in \pr2
meeting at $p_1$, let $p_2$, $p_3$ and $p_4$ be distinct points on $L_1$ away from
$p_1$, let $p_5\ne p_1$ be a point of $L_2$, and let $p_6$ be the point
infinitely near to $p_5$ corresponding to the tangent direction
at $p_5$ along $L_2$. Then $Z=3p_1+2p_2+2p_3+p_4+3p_5+2p_6$ satisfies
the proximity inequalities. Let $e_0,\ldots,e_6$ be the exceptional configuration
coming from the blowing up $X$ of \pr2 at the essentially distinct 
points $p_1,\ldots,p_6$. 

By \ir{eff}(a) any reduced and irreducible curve of 
negative self-intersection is either a component of the total transform 
of a line through any two of the points, a component of the exceptional curve
corresponding to one of the points, or a component of the total transform of
a conic through the points. Thus the classes of reduced and irreducible
curves of negative self-intersection on $X$ are: $e_0-e_1-\cdots-e_4$, $e_0-e_1-e_5-e_6$,
$e_0-e_i-e_5$ for $i\in\{2,3,4\}$, $e_1$, $e_2$, $e_3$, $e_4$, $e_5-e_6$, and $e_6$.
If we let $\C F_d$ denote $de_0-3e_1-2e_2-2e_3-e_4-3e_5-2e_6$,
then for each $d<5$ one can find a sequence
$\C C_1,\ldots,\C C_t$ among the enumerated classes of negative self-intersection
such that each of $\C F_d\cdot \C C_1$, $\ldots$, $(\C F_d-\C C_1-\cdots-\C C_{t-1})\cdot \C C_t$,
and $(\C F_d-\C C_1-\cdots-\C C_t)\cdot e_0$ is negative and hence neither
$\C F_d-\C C_1-\cdots-\C C_t$ nor $\C F_d$ can be the class of an effective divisor.
Thus $h^0(X, \C F_d)=0$ for $d<5$.

In the case that $d=5$, this process of subtracting off the classes of putative fixed components
of negative self-intersection leads to the Zariski decomposition
$\C F_5=\C H+\C N$, where $\C H=2e_0-e_2-e_3-e_5$ is numerically effective
with $h^0(X, \C F_5)=h^0(X,\C H)$, 
and $\C N=3e_0-3e_1-e_2-e_3-e_4-2e_5-2e_6$ with $h^0(X,\C N)=1$.
By \ir{eff}(c), this
process of subtracting off putative components
of negative self-intersection will always either lead to a Zariski decomposition
or to a determination that the original class is not the class of an effective divisor.
By \ir{eff}(b) and \ir{RR}(c),  $h^1(X,\C H)=h^2(X,\C H)=0$,
so $h^0(X,\C F)=h^0(X,\C H)$ is easily computed by Riemann-Roch to be 3.

Similarly, we find: 
that $h^0(X,\C F_d)$ is, respectively, 8, 14, and 23, for $6\le d\le8$;
that in a Zariski decomposition of $\C F_d$ 
the numerically effective part of $\C F_d$ for $6\le d\le 7$ is, respectively,  
$4e_0-e_1-e_2-e_3-2e_5-e_6$,
and $5e_0-e_1-e_2-e_3-2e_5-e_6$; and, adding $e_0$ 
to the numerically effective parts for $5\le d\le 7$
and applying $h^0$, that we get, respectively, 7, 14, and 21.
Thus $\s(\C F_d,e_0)$ is 0 for $d<4$; 3, 1, 0 and 2 for $4\le d\le 7$
(using \ir{cokfact}); and 0 by \ir{myomega} for $d>7$.

Thus, in a resolution $0 \to F_1\to F_0\to I_Z\to 0$ of $I_Z$,
we find that $F_0=R[-5]^3\oplus R[-6]\oplus R[-8]^2$, which allows us as discussed
in \ir{introsect} to determine $F_1=R[-6]^2\oplus R[-7]\oplus R[-9]^2$. \qed

The case of points on a line is a special case of points on a conic
which affords an especially nice and explicit answer, so we present this
case as another example of finding a resolution using our results.

\rem{Example}{linecase}
Let $p_1,\ldots,p_r$ be essentially distinct points
of a line in \pr2; i.e., $e_0-e_1-\cdots-e_r$ is the class of an effective divisor
on the blowing up $X$ of $p_1,\ldots,p_r$. In this case 
finding Zariski decompositions is straightforward, $h^1$ and $h^2$ vanish for
every numerically effective class (by \ir{sg} and \ir{RR}(c)),
and $\C S(\C H,e_0)=0$ for any numerically effective 
class \C H (by \ir{cokvan} or \ir{cokonaconic}), which allows us to determine
a resolution for any $Z=m_1p_1+\cdots+m_rp_r$.

We now make this resolution explicit, leaving details to the reader.
Let $Z=m_1p_1+\cdots+m_rp_r$ be a nontrivial subscheme satisfying
the proximity inequalities. By \ir{reindex}, we may assume that 
$m_1\ge\cdots\ge m_r>0$. From the sequence $m_1\ge\cdots\ge m_r>0$
we get a Young diagram ($r$ columns where the $i$th column is a column of
$m_i$ boxes, $1\le i\le r$), and from this we get the conjugate sequence
$\mu_1\ge \cdots\ge \mu_{m_1}$ ($\mu_i$ is the number of boxes in the $i$th row of the 
Young diagram). We also define $a_i=(i-1)+\mu_i+\cdots+\mu_{m_1}$ for
$1\le i\le m_1$. This gives
$\mu_1+\cdots+\mu_{m_1}=m_1+\cdots+m_r=a_1\ge a_2\ge \cdots\ge a_{m_1}\ge m_1$.
Then the minimal free resolution
of $I_Z$ takes the form
$0\to F_1\to F_0\to I_Z\to 0$, where
$$F_0=R[-a_1]\oplus \cdots \oplus R[-a_{m_1}]\oplus R[-m_1]$$
and
$$F_1=R[-1-a_1]\oplus \cdots \oplus R[-1-a_{m_1}].$$
In particular, we see that a minimal homogeneous generating set for
$I_Z$ has $m_1-m_2+1$ generators of degree $m_1$, and one generator,
of degree $a_i$, for each $1\le i \le m_2$. 

From this resolution we also 
get a particularly nice expression (reformulating and extending
that of [3, Proposition 3.3], which is for {\it distinct\/} points on a line in \pr2)
for the Hilbert function
of $I_Z$, where we follow the convention that $\smallbinom{a}{b}=0$ if $a<b$:
For all $n\ge 0$, 
$h_{I_Z}(n)=\smallbinom{n-m_1+2}{2}+\sum_{1\le i\le m_1}(\smallbinom{n-a_i+2}{2}-
\smallbinom{n-a_i+1}{2})$, or, alternatively, 
$h_{I_Z}(n)=((m_1-m_2+1)\smallbinom{n-m_1+2}{2}-(m_1-m_2)\smallbinom{n-m_1+1}{2})
+\sum_{1\le i\le m_2}(\smallbinom{n-a_i+2}{2}-\smallbinom{n-a_i+1}{2})$.\qed

\irSubsection{Points on a cubic}{cubicsubsect}

We now consider the case of essentially distinct points $p_1,\ldots,p_r$
on a smooth plane cubic $C$; i.e., $p_1$ is a point of $C$,
$p_2$ is a point of the proper transform of $C$, etc.
In this situation,
the monoid EFF of classes of effective divisors on $X$ is controlled by
ker$(\hbox{Cl}(X)\to \hbox{Cl}(D))$ (which we will denote $\Lambda(X,D)$,
or just $\Lambda$ when $X$ and $D$ are clear from context), where $D$ is the proper 
transform to $X$ of $C$ (and hence a section of $-K_X$), and 
$\hbox{Cl}(X)\to \hbox{Cl}(D)$ is the canonical homomorphism
induced by the inclusion $D\subset X$. As is shown in 
[\trans], if one knows $\Lambda$,
then one can determine $h^0(X, \C F)$ for any class \C F on $X$
and one can also effect Zariski decompositions (indeed, one can
determine the free part $(\C F)_m$ of \C F when \C F is the class of an effective divisor).
As discussed in \ir{introsect}, this reduces determining resolutions of ideals
defining fat point subschemes of \pr2 supported at essentially distinct points of
$C$ to determining $\s(\C F,e_0)$ for numerically effective classes \C F.

However, unlike the case of points on a line or conic, there does not seem to 
be a general principle for handling points on a cubic. In particular,
for points on a cubic
$\C S(\C F,e_0)$ need not vanish for every numerically effective class
\C F, even if its linear system of sections is fixed component free or even
base point free. For example, if $p_1,\ldots,p_6$ are distinct general points of \pr2,
then $\C F=5e_0-2e_1-\cdots-2e_6$ is numerically effective and 
its linear system of sections is fixed component and base point free, but 
$h^0(X,e_0)=h^0(X,\C F)=3$ and $h^0(X,\C F+e_0)=10$, 
so $\s(\C F,e_0)\ge 1$. For a more subtle example (one which is not evident
from a simple dimension count), 
consider eight distinct general points $p_1,\ldots,p_8$ of \pr2 and
let $\C F=t(17e_0-6(e_1+\cdots+e_8))$; then, for all $t>0$, \C F is numerically effective and 
its linear system of sections is fixed component and base point free, but 
$\s(\C F,e_0)>0$ (see [\syzconf]). 

Thus we will obtain our results under certain restrictions:
we will restrict either the classes \C F or the surfaces $X$ which we consider. 
More specifically, we will first allow $X$ to be any blowing up of \pr2 at $r\ge 9$
essentially distinct points $p_1,\ldots,p_r$ of a smooth plane cubic,
but only consider {\it uniform classes} (i.e., classes \C F satisfying
$\C F\cdot e_1=\cdots=\C F\cdot e_r$). Second, we will consider
arbitrary classes \C F in the case that $p_1,\ldots,p_r$, $r\ge1$, are
essentially distinct points on a smooth plane cubic $C$ but where $\Lambda$
is as large as possible (i.e., $\Lambda$ is the subgroup $K_X^\perp$
of classes \C F with $\C F\cdot K_X=0$, which in more concrete terms means
$p_1$ is a flex of $C$, and for each $i$,
$p_i$ is the point of the proper transform of $C$ infinitely near to $p_1$).
We begin now by studying uniform classes.

\irSubsubsection{Uniform classes on a blowing up of 
points on a smooth plane cubic}{unifsubsect}

Let $r\ge 9$, and let $p_1,\ldots,p_r$ be essentially distinct
points of a smooth irreducible plane cubic $C$. 
The results in this subsection allow one implicitly to compute
a resolution of the ideal of any fat point subscheme 
of \pr2 of the form $Z=m(p_1+\cdots+p_r)$,
where $m\ge 1$. For this, it is enough to consider uniform classes;
i.e., those of the form
$\C F_n=ne_0-m(e_1+\cdots+e_r)$, where $X$ is the blowing up of
the points $p_1,\ldots,p_r$ and $e_0,\dots,e_r$ is the corresponding
exceptional configuration. Since
$-K_X=3e_0-e_1-\cdots-e_r$, we can write $\C F_n$ as
$te_0-mK_X$, where $t=n-3m$. Of course, our main interest is when
$n$ (and hence $3m+t$) is nonnegative, since otherwise $\C F_n$
is not the class of an effective divisor.

The following Proposition recalls facts from [\trans]
which will be helpful both in our analysis (\ir{cokoncubicA}) of $\s(\C F,e_0)$ for a uniform
class \C F and also in working out complete examples of
resolutions, as in \ir{ristwelve}. 

\prclm{Proposition}{unifFm}{Let $X$ be as in the preceding two paragraphs,
with $\C F=te_0-mK_X$ for some integers $t$ and $m>0$. 
\item{(a)} The class \C F is the class of an effective divisor if and only if $t\ge 0$.
\item{(b)} Say $t\ge0$  and $r=9$.
\itemitem{(i)} If $t>0$, then the linear system of sections of \C F is
fixed component free and $h^1(X, \C F)=0$.
\itemitem{(ii)} If $t=0$, then the linear system of sections of \C F has a fixed
component if and only if $\C F$ is not in $\Lambda$; in particular,
the fixed component free part $(\C F)_m$
of $\C F$ is $-(m-s)K_X$ and we have $h^0(X, (\C F)_m)=\lambda+1$,
where $s$ is the least nonnegative integer
such that $-(m-s)K_X\in \Lambda$ and where $\lambda=0$ if $m=s$ and otherwise
$\lambda=(m-s)/l$, where $l$ is the least positive integer
such that $-lK_X\in \Lambda$.
\item{(c)} Lastly, say $t\ge0$ and $r>9$.
\itemitem{(i)} If $-K_X\cdot \C F>0$, then the linear system of sections of 
\C F is fixed component free and $h^1(X,\C F)=0$.
\itemitem{(ii)} Say $-K_X\cdot \C F\le0$. If $t=0$, then
$(\C F)_m=0$, so say $t>0$ and let
$s$ be the least nonnegative integer such that 
$-K_X\cdot (\C F+sK_X)\ge 0$. If $-K_X\cdot (\C F+sK_X)>0$,
then $(\C F)_m=\C F+sK_X$. If $-K_X\cdot (\C F+sK_X)=0$,
then $(\C F)_m=\C F+sK_X$ and $h^1(X, \C F+sK_X)=1$
if $\C F+sK_X\in\Lambda$, while $(\C F)_m=\C F+(s+1)K_X$
if $\C F+sK_X\not\in\Lambda$.}

\Prf (a) Since $e_0$ and $-K_X$ are the classes of effective divisors, 
\C F is the class of an effective divisor if $t\ge 0$.
Conversely, note that $-K_X$ is the class of an irreducible curve 
of self-intersection $9-r$, and that
$-K_X\cdot \C F=3t+(9-r)m$, which is negative if $t<0$ and $r\ge9$.
Now, $-K_X$ is numerically effective for $r=9$, but, if $t<0$, meets \C F negatively,
which implies that \C F cannot be the class of an effective divisor.
For $r>9$ and $t<0$, if \C F is the class of an effective divisor, then $-K_X$ is a fixed component so
$\C F+K_X=te_0-(m-1)K_X$ is the class of an effective divisor. Iterating we eventually obtain
the contradiction that $te_0$ is the class of an effective divisor, which is absurd if $t<0$.

Items (b,c)(i) follow by Theorems 1.1(b) and 3.1 of [\trans].
Items (b,c)(ii) follow by Proposition 1.2 and Theorem 3.1 of [\trans]. \qed

We now compute $\s(\C F_n,e_0)$ for each $n>0$. 
By \ir{unifFm}, the fixed component free part $(\C F_n)_m$ of $\C F_n$ is
a uniform class, hence so is $e_0+(\C F_n)_m$. 
Thus, using \ir{unifFm}, it is enough by \ir{cokfact} to compute
$\s((\C F_n)_m,e_0)$, for which the following theorem suffices.

\prclm{Theorem}{cokoncubicA}{Let $\C F=(\C G)_m$ for some uniform class \C G of 
an effective divisor, where $X$ is the blowing up of $r\ge 9$ essentially 
distinct points of a smooth plane cubic.
\item{(a)} If $-K_X\cdot \C F>1$, then $\C S(\C F,e_0)=0$.
\item{(b)} If $-K_X\cdot \C F=1$, then $\s(\C F,e_0)=1$.
\item{(c)} If $-K_X\cdot \C F=0$, then $\C S(\C F,e_0)=0$ unless either:
$r=10$, in which case $\s(\C F,e_0)=1$;
or $r=9$, in which case $\C F=-abK_X$ and $\s(\C F,e_0)=3b(a-1)$,
where $a=b=0$ if $\C F=0$, and otherwise
$a$ is the least positive integer such that $-aK_X\in \Lambda$
and $b$ is some nonnegative integer.}

\Prf Write \C F as $te_0-mK_X$ for some nonnegative integers $t$ and $m$.

(a) Note that $-K_X\cdot(te_0-mK_X)>1$
implies $-K_X\cdot(te_0-sK_X)>1$
for all $0\le s\le m$. We will induct on $s$, starting with
the obvious fact that $\C S(te_0,e_0)=0$ for $t\ge 0$. So now we may assume
that $0<s\le m$, and that $\C S(te_0-(s-1)K_X,e_0)=0$. 
By \ir{unifFm}, $h^1(X,te_0-(s-1)K_X)=h^1(X,te_0-(s-1)K_X+e_0)=0$, and
$h^1(X,e_0+K_X)=h^1(X,-e_0)=h^1(\pr2,-e_0)=0$, so 
by \ir{Mumford}(a,b) we have the exact sequence
$\C S(te_0-(s-1)K_X,e_0)\to\C S(te_0-sK_X,e_0)\to\C S((te_0-sK_X)\otimes\C O_D,e_0\otimes\C O_D)\to 0$,
where the leftmost term vanishes by induction, and $\C S((te_0-sK_X)\otimes\C O_D,e_0\otimes\C O_D)=0$
by \ir{Mumford}(c), since $(te_0-sK_X)\otimes\C O_D$ has degree at least $2g$
and $e_0\otimes\C O_D$ has degree $2g+1$, where $g=1$ is the genus of $D$.
Thus $\C S(te_0-sK_X,e_0)=0$ follows by exactness.

(b) Under the given hypotheses, we must have
$m>0$ and $K_X^2<0$, hence 
$-K_X\cdot (\C F+K_X)>1$ and $\C S(\C F+K_X,e_0)=0$ by (a). 
As in (a), we have an exact sequence
$\C S(\C F+K_X,e_0)\to\C S(\C F,e_0)\to\C S(\C F\otimes\C O_D,e_0\otimes\C O_D)\to0$.
Applying it gives 
$\s(\C F,e_0) =\s(\C F\otimes\C O_D,e_0\otimes\C O_D)$,
but $h^0(D, \C F\otimes\C O_D)=1$ since $\C F\otimes\C O_D$ has degree 1.
Thus $H^0(D,\C F\otimes\C O_D)\otimes H^0(D, e_0\otimes\C O_D)\to H^0(D, (\C F+e_0)\otimes\C O_D)$
is injective and one easily computes $\s(\C F\otimes\C O_D,e_0\otimes\C O_D)=1$.

(c) If $\C F=0$, then $\s(\C F,e_0)=0$ is clear, so assume $\C F\ne0$.
Since $-K_X\cdot\C F=0$, we see $m>0$ and $\C F+K_X$ is
numerically effective. Since $\C F=(\C F)_m$ with $-K_X\cdot \C F=0$,
we see $\C F\in \Lambda$, so $\C F\otimes\C O_D=\C O_D$ and $h^0(D,\C O_D\otimes\C F)=1$,
and as above we see $h^1(X,e_0+K_X)=h^1(X,\C F+K_X+e_0)=0$.
If $-K_X^2>0$, then also $h^1(X,\C F+K_X)=0$ by \ir{unifFm}(c). 
If, however, $-K_X^2=0$, then using \ir{unifFm}(b) we see
$h^0(X,\C F+K_X)+h^0(D,\C F\otimes\C O_D)=h^0(X,\C F)$, so
$0\to \C F+K_X\to \C F\to\C F\otimes\C O_D\to 0$ is exact 
on global sections. In either case, then, by \ir{Mumford}(b,c) we have an exact sequence
$\C R(\C F\otimes\C O_D,e_0)\to 
\C S(\C F+K_X,e_0)\to\C S(\C F,e_0)\to\C S(\C F\otimes\C O_D,e_0\otimes\C O_D)\to0$.

Since $\C F\otimes\C O_D=\C O_D$, we see
$\C R(\C F\otimes\C O_D,e_0)=0$ and $\s(\C F\otimes\C O_D,e_0\otimes\C O_D)=0$.
Thus we have $\s(\C F,e_0)=\s(\C F+K_X,e_0)$; now $\s(\C F+K_X,e_0)$ equals 0
by (a) if $r>10$, and 1 by (b) if $r=10$. 
If $r=9$, then $\C F=-mK_X$. By \ir{unifFm}(b), $m=ab$, where $a$ is the
least positive integer such that $-aK_X\in\Lambda$. Likewise, 
$(\C F+K_X)_m=-a(b-1)K_X$ and $h^0(X, \C F+K_X+e_0)-h^0(X, (\C F+K_X)_m+e_0)=3a-3$;
by induction, we may assume $\s(\C F+K_X,e_0)=3(b-1)(a-1)$, so 
$\s(\C F+K_X,e_0)=3b(a-1)$ by \ir{cokfact}. \qed

\rem{Example}{ristwelve} We now give an example using our results
to obtain an explicit resolution. Consider subschemes $Z\subset\pr2$ of the form 
$Z=m(p_1+\cdots+p_{12})$, where $p_1,\ldots,p_{12}$ are essentially distinct 
points on a smooth cubic. Let $X$ be the blowing up of the points,
let $e_0,\ldots,e_{12}$ be the associated exceptional configuration, let $D$ be the proper
transform of the cubic, and assume that the kernel $\Lambda$ of $\hbox{Cl}(X)\to \hbox{Cl}(D)$
is trivial. Then the minimal free resolution of $I_Z$ is:
$$0\to \mylim{$\bigoplus$}{$\scriptstyle1\le i\le m$}{5}R[-3m-i-2]^3\to 
R[-3m]\oplus(\mylim{$\bigoplus$}{$\scriptstyle1\le i\le m$}{5}R[-3m-i-1]^3)\to I_Z\to 0.$$

We sketch a proof. Let $\C F_n$ denote $ne_0-m(e_1+\cdots+e_{12})$. By \ir{unifFm},
$h^0(X,\C F_n)=0$ for $n<3m$, and for $n\ge 3m$
we see that $(\C F_n)_m$ is regular, equal to $0$ if $n=3m$
and otherwise to $(n-3m-1)(-K_X)+(n-3m)e_0$.
Now, by \ir{cokfact} and \ir{cokoncubicA}, we can compute that $\C S(\C F_n,e_0)$
vanishes for $n\le 3m-2$ or $n\ge 4m+1$, equals 1 for $n=3m-1$,
0 for $n=3m$ and 3 for $3m+1\le n\le 4m$.
Thus the first syzygy module $F_0$ in the resolution is as claimed;
now we merely need to check that the difference of the Hilbert functions of
$F_0$ and $I_Z$ coincides with the Hilbert function of 
$\bigoplus_{1\le i\le m}R[-(3m+i+2)]^3$.    \qed

\irSubsubsection{Points infinitely near a flex}{flexsubsect}

Let $p_1,\ldots,p_r$ be essentially distinct points on a smooth plane cubic $C$,
let $X$ be the blowing up of the points with $e_0,\ldots,e_r$ the corresponding exceptional
configuration, and let $D$ denote the proper transform of $C$ on $X$.
In this subsection we shall always assume that
$p_1$ is a flex of $C$, and $p_i$ for $i>1$ is infinitely near to $p_1$. This
is equivalent to $r_0=e_0-e_1-e_2-e_3$ and $r_i=e_i-e_{i+1}$, $i>0$, all being in the kernel $\Lambda$ of
$\hbox{Cl}(X)\to \hbox{Cl}(D)$. But $r_0$ and $r_i$, $i>0$,
give a basis for the subspace $K^\perp$ of $\hbox{Cl}(X)$ of classes perpendicular to $K_X$
and so the requirements that $p_1$ be a flex of $C$, and $p_i$ for $i>1$ be points of proper transforms
of $C$ and infinitely near to $p_1$ are equivalent to $\Lambda=K^\perp$. Thus we may equivalently
say that $p_1,\ldots,p_r$ are essentially distinct points on a smooth plane cubic such that
$\Lambda=K^\perp$. 

Zariski decompositions and $h^0(X, \C F)$ for any $\C F\in\hbox{Cl}(X)$
can be computed using [\trans].
Thus, as discussed in \ir{introsect}, to compute a resolution for $I_Z$
for any fat point subscheme $Z$ supported at $p_1,\ldots,p_r$ 
it suffices to determine $\s(\C H, e_0)$
for numerically effective classes \C H, which we do below. 
Thus the results below allow one to work out
a resolution for the ideal of any fat point subscheme
supported at $p_1,\ldots,p_r$. 

The values of $\s(\C F, e_0)$ that we obtain are related to the structure of NEFF, so
we begin by describing the cone NEFF
of numerically effective divisors, for which we use the following notation:
$\C H_0=e_0$, $\C H_1=e_0-e_1$, $\C H_2=2e_0-e_1-e_2$, and $\C H_i=3e_0-e_1-\cdots-e_i$,
$i\ge 3$. In addition, we will need to recall facts about $h^1$ for numerically effective
classes, which will also be helpful to readers 
moved to work out examples of resolutions of fat point ideals for fat point subschemes
supported at $p_1,\ldots,p_r$.

\prclm{Lemma}{neffkerall}{Let \C F be a class on $X$; then:
\item{(a)} $\C F\in\hbox{NEFF}$ if and only if \C F 
is a nonnegative (integer) linear combination 
of $\C H_i$, $i\ge 0$, such that $-\C F\cdot K_X\ge 0$; and 
\item{(b)} for $\C F\in\hbox{NEFF}$ we have
$$h^1(X,\C F)=\cases{0, &if $-K_X\cdot\C F\ge 1$;\cr
                                1, &if $\C F^2>-K_X\cdot\C F=0$;\cr
                    \C F\cdot e_1, &if $\C F^2=-K_X\cdot\C F=0$.\cr}$$
\noindent Moreover, writing $\C F\in\hbox{NEFF}$ as $\C F=\sum_ia_i\C H_i$,
the linear system of sections of \C F has a nonempty base locus if and only if:
$\C F=\C H_8$ or $-K_X\cdot\C F=1$ but $\C F\ne \C H_8+a_9\C H_9$ (in this case the base 
locus is a single point if $j=r$, and it is the divisor $E_{j+1}$ if $j<r$, 
where $j$ is the greatest subscript with $a_j>0$); 
or $\C F=\C H_8+a_9\C H_9$ with $a_9>0$ (in which case $E_9$ is the base locus);
or $\C F=\C H_8+a_9\C H_9+\C H_{10}$ (in which case $E_9-E_{10}$
is the base locus).}

\Prf (a) See [\trans].  

(b) See [\trans]. We remark that $\C F^2=-K_X\cdot\C F=0$ only occurs
for $\C F=s\C H_9$, for some $s\ge 0$, in which case $s=\C F\cdot e_1$.  \qed

By \ir{neffkerall}, NEFF is contained in the nonnegative subsemigroup of 
$\hbox{Cl}(X)$ generated by $\C H_0,\ldots,\C H_r$. It turns out to be convenient
to distinguish two types. Let $\C F=\sum_ia_i\C H_i$ be a numerically effective class
(hence $a_i\ge 0$ for all $i$). 
We say \C F is of {\it type I\/} if $a_i>0$ for some $i<8$; otherwise, we say \C F is of
{\it type II\/} (i.e., either $a_i=0$ for $i\ne 9$ or the least index $l$ such that
$a_l\ne 0$ is 8).
We first consider classes of type I.

\prclm{Theorem}{cokoncubicB}{Let \C F
be a type I numerically effective class on $X$,
and, writing \C F as a nonnegative linear combination $\C F=\sum_ia_i\C H_i$,
let $j$ be the largest index with $a_j>0$. Then 
$\s(\C F,e_0)=0$ unless either $-K_X\cdot\C F=1$, or $-K_X\cdot\C F=0$ and $j=10$,
in which cases $\s(\C F,e_0)=1$.}

\Prf Clearly, $\C S(\C O_{X_l},e_0)=0$, so we may assume that $\C F\ne 0$.
Let $l$ be the least index $i$ such that $a_i>0$; then $l<8$.
By induction, we may assume our result is true for $\C F-\C H_j$. 
For each $i$, let $X_i$ be the blowing up of $p_1,\ldots,p_i$ and let
$D_i$ be the proper transform of $C$ to $X_i$.
By duality and the fact that $K_{X_i}+\C H_i$ is always a multiple
of $e_0$, we see $h^1(X_i,e_0-\C H_i)=0$, and, applying \ir{neffkerall},
we see $h^1(X_j,\C F-\C H_j)$ and $h^1(X_j,\C F-\C H_j+e_0)$ vanish. Thus 
(suppressing the subscripts on $X_j$ and $D_j$) we have the usual exact sequence
$\C R(\C F\otimes\C O_D,e_0)\to\C S(\C F-\C H_j,e_0)\to\C S(\C F,e_0)\to
\C S(\C F\otimes\C O_D,e_0\otimes\C O_D)\to0$.

If $-K_X\cdot \C F\ge2$, then $\C S(\C F\otimes\C O_D,e_0\otimes\C O_D)=0$ 
by \ir{Mumford}(c), and either $\C F-\C H_j=0$, or
$-K_X\cdot (\C F-\C H_j)\ge -K_X\cdot \C H_l\ge2$, so 
$\C S(\C F-\C H_j,e_0)=0$ by induction, and, applying the exact sequence,
$\C S(\C F,e_0)=0$.     

If $-K_X\cdot \C F<2$, then $j>9$ and, as in the proof of \ir{cokoncubicA}(b,c),
$\C R(\C F\otimes\C O_D,e_0)=0$. Now, $\s(\C F\otimes\C O_D,e_0\otimes\C O_D)$ is 0 if
$-K_X\cdot \C F=0$ and 1 if $-K_X\cdot \C F=1$, so, 
by exactness,  $\s(\C F,e_0)=\s(\C F-\C H_j,e_0)-K_X\cdot\C F$
if $0\le -K_X\cdot\C F\le 1$, and applying the inductive hypothesis to
$\C F-\C H_j$ gives the result.
\qed 

We now consider numerically effective classes \C F of type II. Note that if
$\C F=\C H_8+b_9\C H_9+b_{10}\C H_{10}$ with either $b_9$ or $b_{10}$
positive, then the linear system of sections of \C F has a fixed component,
making this a case we do not need to consider (because then $\s(\C F,e_0)=b_9+1$
follows from \ir{cokfact} and the following theorem). Thus the only type II
classes we need to consider are $\C H_8$, $b_9\C H_9$ and numerically effective
classes $\sum_{i\ge 8}b_i\C H_i$ with $b_8>1$.

\prclm{Theorem}{cokoncubicC}{Let $\C F=\sum_{i\ge 8}b_i\C H_i$ be a numerically 
effective class on $X$ (hence $b_i\ge 0$ for each $i$), and, if $\C F\ne 0$, let 
$j$ be the greatest index $i$ such that $b_i>0$.
\item{(a)} We have $\s(\C H_8,e_0)=1$ and $\s(b_9\C H_9,e_0)=0$.
\item{(b)} If $b_8>1$, then $\s(\C F,e_0)=1$ unless either
$-K_X\cdot\C F=1$, or $-K_X\cdot\C F=0$ and $j=10$,
in which cases $\s(\C F,e_0)=2$.}

\Prf As in the proof of \ir{cokoncubicB} (and again 
suppressing subscripts on $X_j$ and $D_j$),
\ir{Mumford} gives an exact sequence
$\C R(\C F\otimes\C O_D,e_0)\to\C S(\C F-\C H_j,e_0)\to\C S(\C F,e_0)\to
\C S(\C F\otimes\C O_D,e_0\otimes\C O_D)\to0$.

(a) Since $\C F=b_9\C H_9$ is uniform (on the blowing up
of \pr2 at $p_1,\ldots,p_9$), we have $\C S(\C F,e_0)=0$ by 
\ir{cokoncubicA}(c), as required, so assume  
$\C F=\C H_8$. Clearly, $\C S(\C O_X,e_0)=0$, 
so, by the exact sequence with $j=8$, $\s(\C H_8,e_0)=\s(\C H_8\otimes\C O_D,e_0\otimes\C O_D)$,
and, as in the proof of \ir{cokoncubicA}(b), 
$\s(\C H_8\otimes\C O_D,e_0\otimes\C O_D)=1$, so $\s(\C H_8,e_0)=1$,
as required. 

(b) We first show the $j=8$ case implies the others. So suppose $j=8$ and
assume $\s(b\C H_8,e_0)=1$ for all $b\ge0$. 
Since  $\C S((b+1)\C H_8\otimes\C O_D,e_0\otimes\C O_D)=0$ for all $b\ge 1$
by \ir{Mumford}(c), by the exact sequence above
the homomorphism $\C S(b\C H_8,e_0)\to \C S((b+1)\C H_8,e_0)$
is an isomorphism for all $b>0$.

Now suppose $j>8$.
By induction, we have an isomorphism $\C S(b_8\C H_8,e_0)\to \C S((\sum_{i\ge 8}b_i)\C H_8,e_0)$,
and it factors as $\C S(b_8\C H_8,e_0)\to \C S(\C F,e_0)\to \C S((\sum_{i\ge 8}b_i)\C H_8,e_0)$,
which implies $\s(\C F,e_0)\ge 1$. 
If $-K_X\cdot\C F>1$, then $\s(\C F\otimes\C O_D,e_0\otimes\C O_D)=0$
by \ir{Mumford}(c), so $\C S(\C F-\C H_j,e_0)\to\C S(\C F,e_0)$ is surjective.
By induction we may assume that $\s(\C F-\C H_j,e_0)=1$, so $\s(\C F,e_0)\le 1$
which with $\s(\C F,e_0)\ge 1$ gives $\s(\C F,e_0)=1$.
If, however, $-K_X\cdot\C F=1$, then $b_8>1$ implies $j>9$, so 
$-K_X\cdot (\C F-\C H_j)>1$, and by induction we may assume
$\s(\C F-\C H_j,e_0)=1$. But, as above,
the isomorphism $\C S(b_8\C H_8,e_0)\to\C S((\sum_{8\le i}b_i)\C H_8,e_0)$
factors through $\C S(\C F-\C H_j,e_0)\to\C S(\C F,e_0)$, so
the latter is injective. The exact sequence with
$\s(\C F\otimes\C O_D,e_0\otimes\C O_D)=1$ now gives $\s(\C F,e_0)=2$.
For the case that $-K_X\cdot\C F=0$, we have $\C F\otimes\C O_D=\C O_D$ so
$\C R(\C F\otimes\C O_D,e_0)=0=\C S(\C F\otimes\C O_D,e_0\otimes\C O_D)$.
Thus  $\s(\C F-\C H_j,e_0)=\s(\C F,e_0)$;
but $j\ge 10$ so by what we have already done we see
$\s(\C F-\C H_j,e_0)=1$ if $j>10$ and
$\s(\C F-\C H_j,e_0)=2$ if $j=10$. This finishes the
proof, modulo showing $\s(b\C H_8,e_0)=1$ for all $b>0$,
to which we now proceed.

We already have observed that $\s(\C H_8,e_0)=1$ so we may assume $b\ge 2$.
The pencil of cubics through $p_1,\ldots,p_8$ includes our smooth cubic $C$ and
a triple line $3L$ tangent to $C$ at the flex $p_1$. 
If $C$ is not supersingular in characteristic 3, then 
$\hbox{Pic}(C)$ has nontrivial 3-torsion so there is a second flex. We will denote the
line on $X$ tangent to $D$ at this second flex by $S$, and we will
choose projective coordinates $x,y,z$ on \pr2 such that $z=0$ defines $L$,
$y=0$ defines the line through $p_1$ and the second flex, and
$x$ defines the line tangent to $C$ at the second flex (whose transform on $X$ is thus $S$). 
By scaling $x$, $y$ and $z$, the homogeneous form $f$ defining $C$ 
can be taken to be $y^3+x^2z+xz^2+axyz$, for some $a\in k$. 
If $C$ is supersingular in characteristic 3, then we can take $f$ to be 
$y^3+x^2z+yz^2$ (see Proposition IV.4.21 and Example IV.4.23.1 of [14]).

By \ir{Mumford}(a), we have the exact sequence
$\C S(\C O_X,b\C H_8)\to\C S(e_0,b\C H_8)\to\C S(e_0\otimes\C O_S,b\C H_8)\to0$.
Clearly $\C S(\C O_X,b\C H_8)=0$,
so $\s(e_0,b\C H_8)=\s(e_0\otimes\C O_S,b\C H_8)$. Taking $y$ and $z$ for projective coordinates
on $S$, $H^0(S,e_0+b\C H_8)$ can be identified with the vector space with basis consisting of 
the monomials in $y$ and $z$ of degree $3b+1$. By
explicitly computing the image of $H^0(S,e_0)\otimes H^0(X,b\C H_8)\to H^0(S,e_0+b\C H_8)$,
we will see that $\s(e_0\otimes\C O_S,b\C H_8)=1$.

Under the identification of $H^0(X, de_0)$ with $H^0(\pr2,de_0)$, which we regard
as the space of monomials in $x$, $y$ and $z$ of degree $d$, we can 
for any class \C G on $X$ regard $H^0(X,\C G)$ as a subspace of 
a space $H^0(\pr2,(e_0\cdot\C G)e_0)$ of monomials of degree $e_0\cdot\C G$.
From this point of view, 
$\{f,z^3\}$ is a basis for $H^0(X,\C H_8)$. (To say this differently, 
$f$ and $z^3$ are forms  passing through $p_1,\ldots,p_8$ with multiplicity at least 1
at each point, so each is an element of $H^0(\pr2,3e_0)$ corresponding
to an element of $H^0(X,\C H_8)$. Since $f$ and $z^3$ are linearly independent
in $H^0(\pr2,3e_0)$ and $h^0(X,\C H_8)=2$, they correspond to a basis. Of course,
that $f$ passes through each of $p_1,\ldots,p_8$ with multiplicity at least 1
is obvious; to see the same for $z^3$, note that $\C H_8=3r_0+2r_1+4r_2+6r_3+5r_4+4r_5+3r_6+2r_7+e_8$
and that each of $r_1,\ldots,r_8$ and $e_8$ is the class of a component of $E_1$, while
$r_0$ is the class of the proper transform of $z=0$.) Similarly,
$\{f^2,fz^3,z^6, yz^5\}$ is a basis for $H^0(X,2\C H_8)$
(to see that $yz^5\in H^0(\pr2,6e_0)$ corresponds to
a section of $2\C H_8$, note that $2\C H_8=(e_0-e_1)+5r_0+4r_1+7r_2+10r_3+8r_4+6r_5+4r_6+2r_7$,
so the $y$ comes as a section of $e_0-e_1$ and the $z^5$ from $5r_0$), 
$\{f^3, f^2z^3, fz^6, fyz^5, z^9, yz^8, xz^8\}$ is a basis for $H^0(X,3\C H_8)$ 
(to see that $xz^8\in H^0(\pr2,9e_0)$ corresponds to
a section of $3\C H_8$, note that $3\C H_8=e_0+8r_0+5r_1+10r_2+15r_3+12r_4+9r_5+6r_6+3r_7$,
so the $x$ comes as a section of $e_0$ and the $z^8$ from $8r_0$), and 
$\{f^4, f^3z^3, f^2z^6, f^2yz^5, fz^9, fyz^8, fxz^8, z^{12}, yz^{11}, xz^{11}, y^2z^{10}\}$ 
is a basis for $H^0(X,4\C H_8)$.

Also note that $\C S(2\C H_8,3\C H_8)=0$ (one checks that
$fH^0(X,4\C H_8)\subset H^0(X,5\C H_8)$ is in the image of 
$H^0(X,2\C H_8)\otimes H^0(X,3\C H_8)$ explicitly and that 
$\C S(2\C H_8\otimes \C O_D,3\C H_8\otimes \C O_D)=0$), so 
$H^0(X,5\C H_8)$ is spanned by products of sections
of $H^0(X,2\C H_8)$ and $H^0(X,3\C H_8)$.
One can now check explicitly for $1\le b\le 5$ that the image of the restriction homomorphism
$H^0(X,b\C H_8)\to H^0(S,b\C H_8)$ is contained in
the space spanned by $\{z^{3b}, yz^{3b-1},\ldots, y^{3b-3}z^3, (y^3+\alpha yz^2)^b\}$,
where $\alpha$ is 1 if $C$ is supersingular in characteristic 3, and 0 otherwise. 
By Proposition III.1 of [\bir], $b\C H_8$ is normally generated for $b\ge 3$,
so $\C S((b-3)\C H_8,3\C H_8)=0$ for all $b>5$. Thus, for $b\ge 6$,
$H^0(X,b\C H_8)$ is spanned by products of forms from  
$H^0(X,3\C H_8)$ and $H^0(X,(b-3)\C H_8)$, from which it is easy to
conclude that the image of the restriction homomorphism
$H^0(X,b\C H_8)\to H^0(S,b\C H_8)$ is contained for all $b$ in
the space spanned by $\{z^{3b}, yz^{3b-1},\ldots, y^{3b-3}z^3, (y^3+\alpha yz^2)^b\}$. 
But an easy argument shows that the subspace of $H^0(S, e_0+b\C H_8)$ spanned by products
of elements from $\{z^{3b}, yz^{3b-1},\ldots, y^{3b-3}z^3, (y^3+\alpha yz^2)^b\}$
and $\{y,z\}$ never includes $y^{3b-2}z^2$; i.e., $\s(e_0\otimes\C O_S, b\C H_8)\ge 1$.
To see that $\s(e_0\otimes\C O_S, b\C H_8)=1$, recall that we have an exact sequence
$\C S((b-1)\C H_8,e_0)\to \C S(b\C H_8,e_0)\to 
\C S(b\C H_8\otimes\C O_D,e_0\otimes\C O_D)\to 0$
with $\C S(b\C H_8\otimes\C O_D,e_0\otimes\C O_D)=0$ for $b\ge 2$
by \ir{Mumford}(c). Since $\s(\C H_8,e_0)=1$ by (a),
we have by induction that $\s(e_0\otimes\C O_S, b\C H_8)=
\s(b\C H_8,e_0)\le 1$ for all $b\ge 1$. \qed

\vskip-\baselineskip
\vskip-\baselineskip
\References

\bibitem{1} Artin, M. {\it Some numerical criteria for contractability of curves
on algebraic surfaces}, 
Amer.\ J.\  Math.\ 84 (1962), 485--497.

\bibitem{2} Catalisano, M.\ V. {\it ``Fat'' points on a conic}, 
Comm.\ Alg.\  19(8) (1991), 2153--2168.

\bibitem{3} Davis, E.\ D.\ and Geramita, A.\ V. {\it The Hilbert
function of a special class of 1-dimensional Cohen-Macaulay
graded algebras},
Queen's papers in Pure and Applied Mathematics, no. 67,
The Curves Seminar at Queen's, vol. III (1984).

\bibitem{4} Davis, E.\ D., Geramita, A.\ V., and Maroscia, P. {\it Perfect
Homogeneous Ideals: Dubreil's Theorems Revisited},
Bull.\ Sc.\ math., $2^e$ s\'erie, 108 (1984), 143--185.

\bibitem{5} Gimigliano, A. {\it Our thin knowledge of fat points},
Queen's papers in Pure and Applied Mathematics, no. 83,
The Curves Seminar at Queen's, vol. VI (1989).

\bibitem{\trans} Harbourne, B. {\it Complete linear systems on rational surfaces}, 
Trans.\ A.\ M.\ S.\ 289 (1985), 213--226. 

\bibitem{\duke} Harbourne, B. {\it Blowings-up of \pr 2 and their blowings-down}, 
Duke Math. J. 52 (1985), 129--148.

\bibitem{\vanc} Harbourne, B. {\it The geometry of rational surfaces and Hilbert
functions of points in the plane},
Can.\ Math.\ Soc.\ Conf.\ Proc.\ 6 
(1986), 95--111.

\bibitem{\ravello} Harbourne, B. {\it Points in Good Position in \pr 2}, in:
Zero-dimensional schemes, Proceedings of the
International Conference held in Ravello, Italy, June 8--13, 1992,
De Gruyter, 1994.

\bibitem{\mtnwest} Harbourne, B. {\it Rational surfaces with $K^2>0$}, preprint, 1994, to appear.

\bibitem{\ars} Harbourne, B. {\it Anticanonical rational surfaces}, preprint (available via my Web
page), 1994.

\bibitem{\bir} Harbourne, B. {\it Birational morphisms of rational surfaces}, preprint, 1994.

\bibitem{\syzconf} Harbourne, B. {\it Generators for Symbolic Powers of 
Ideals Defining General Points of \pr2}, preprint (available via my Web
page), 1995.

\bibitem{14} Hartshorne, R. Algebraic Geometry. Springer-Verlag, 1977.

\bibitem{15} Hirschowitz , A.
{\it Une conjecture pour la cohomologie 
des diviseurs sur les surfaces rationelles generiques},
Journ.\ Reine Angew.\ Math. 397
(1989), 208--213.

\bibitem{16} Iarrobino, A. {\it Inverse system of a symbolic power, III: thin algebras
and fat points}, preprint 1994.

\bibitem{17} Lorenzini, A. {\it The minimal resolution conjecture}, J.\ Alg.
156 (1991), 5--35.

\bibitem{18} Mumford, D. {\it Varieties defined by quadratric equations},
in: Questions on algebraic varieties, Corso C.I.M.E. 1969 Rome: Cremonese,
1969, 30--100.

\bibitem{19} Ramanujam, C.\ P. {\it Supplement to the article 
	``Remarks on the Kodaira vanishing theorem''},
Jour.\ Ind.\ Math.\ Soc. 38 (1974), 121--124.

\vfil\eject

\font\bigrm=cmr12 at 14.4pt
\font\bigi=cmmi12 at 14.4pt
\font\bigscriptrm=cmr12 at 12pt
\font\bigscripti=cmmi12 at 12pt
\font\bigsy=cmsy10 at 14.4pt
\font\bigbf=cmbx12 at 14.4pt
\font\bigtt=cmtt12 at 14.4pt
\font\bigit=cmti12 at 14.4pt
\font\bigsl=cmsl12 at 14.4pt

\def\bigfont{\def\rm{\fam0\bigrm}
\textfont0=\bigrm \scriptfont0=\bigscriptrm \scriptscriptfont0=\tenrm
\textfont1=\bigi \scriptfont1=\bigscripti \scriptscriptfont1=\teni
\def\mit{\fam1} \def\oldstyle{\fam1\teni}
\textfont2=\bigsy 
\def\cal{\fam2}
\def\it{\fam\itfam\bigit} 
\textfont\itfam=\bigit
\def\sl{\fam\slfam\bigsl} 
\textfont\slfam=\bigsl
\def\bf{\fam\bffam\bigbf} 
\textfont\bffam=\bigbf 
\def\tt{\fam\ttfam\bigtt} 
\textfont\ttfam=\bigtt
\normalbaselineskip=17.28pt
\setbox\strutbox=\hbox{\vrule height12.25pt depth5pt  width0pt}%
\normalbaselines\rm}

{\bigfont
\noindent\hbox to\hsize{Addendum to:  \hfil December 5, 2003}
\address
Free Resolutions of Fat Point Ideals on ${\bf P}^2$
Brian Harbourne
J. Pure Appl. Alg. 125, 213--234 (1998)
\endaddress
}
\vskip\baselineskip

This short note adds some details to the proof of \ir{keylemma} which
were not included in the proof in the originally published version of this paper.
The changes are indicated by left and right indentation.
\vskip\baselineskip

\prclm{Lemma}{keycor}{Let $X$ be a smooth projective 
rational surface
and let \C N be the class of a nontrivial effective divisor $N$ 
on $X$. If $\C N+K_X$ is
not the class of an effective divisor and \C F meets every component of \C N nonnegatively, 
then $h^0(N, \C F)>0$, $h^1(N, \C F)=0$, $N^2+N\cdot K_X<-1$,
and every component $M$ of $N$ is a smooth rational curve (of negative
self-intersection, if $M$ does not move).}

\prclm{Lemma}{keylemma}{Let $X$ be a smooth projective 
rational surface, and let \C N be the class of an effective divisor $N$ on $X$
such that $h^0(X, \C N+K_X)=0$. If \C F and \C G are the restrictions to $N$
of divisor classes $\C F'$ and $\C G'$ on $X$ which meet each component of $N$ nonnegatively,
then $\C S(\C F,\C G)=0$.}

\Prf To prove the lemma, induct on the sum $n$ of 
the multiplicities of the components of $N$.
By \ir{keycor}, $h^1(N, \C O)=0$ and every component
of $N$ is a smooth rational curve. Thus the case $n=1$
is trivial (since then $N=\pr 1$, and
the space of polynomials of degree $f$ in two variables 
tensor the space of polynomials of degree $g$ in two variables 
maps onto the space of polynomials of degree $f+g$). So say $n>1$.

As in the proof
of Theorem 1.7 of [1] (or see the proof of
Lemma II.6 of [\mtnwest]), $N$ has a component $C$ such that 
$(N-C)\cdot C\le 1$. Let $L$ be the effective divisor $N-C$
and let \C L be its class.
Thus we have an exact sequence 
$0\to \C O_C\otimes(-\C L)\to \C O_N\to \C O_L\to 0$.
\vskip\baselineskip

{\narrower\smallskip\parindent=0in
To see this, apply the snake lemma to 
$$\matrix{
0 & \to & \C O_X(-N) & \to & \C O_X & \to & \C O_N & \to & 0\cr
  &     & \downarrow &     & \downarrow & & \downarrow & & \cr
0 & \to & \C O_X(-L) & \to & \C O_X & \to & \C O_L & \to & 0\cr}
$$
to see that the kernel of $\C O_N \to \C O_L$ is just
the cokernel of $\C O_X(-N) \to \C O_X(-L)$, which is just
$\C O_C\otimes \C O_X(-L)$, which we may write as $\C O_C(-L)$.
\smallskip}
\vskip\baselineskip

Now, $-L\cdot C\ge -1$, and both 
$\C F'$ and $\C G'$ meet $C$ nonnegatively. We may assume
$\C F'\cdot C\ge\C G'\cdot C$, otherwise reverse the roles of 
$\C F'$ and $\C G'$. Since
$C=\pr 1$, we see that $h^1(C, \C O_C\otimes(\C F'-\C L))$,
$h^1(C, \C O_C\otimes(\C G'-\C L))$ and 
$h^1(C, \C O_C\otimes(\C F'+\C G'-\C L))$ all vanish.
An argument similar to that used to prove \ir{Mumford}(a, b) 
now shows that we have an exact sequence
$\C S(\C O_C\otimes(\C F'-\C L),\C O_C\otimes\C G')\to
\C S(\C F,\C G)\to\C S(\C O_L\otimes\C F,\C O_L\otimes\C G)\to0$.
\vskip\baselineskip

{\narrower\smallskip\parindent=0in
What is actually clear here is that we have
$\C S(\C O_C\otimes(\C F'-\C L),\C G)\to
\C S(\C F,\C G)\to\C S(\C O_L\otimes\C F,\C G)\to0$.
Since $h^1(C, \C O_C\otimes(\C G'-\C L))=0$, we know
$\C G\to \C O_L\otimes\C G$ is surjective on global sections,
and hence that $\C S(\C O_L\otimes\C F,\C G)$ is the same as
$\C S(\C O_L\otimes\C F,\C O_L\otimes\C G)$.
What needs additional justification here is that
$\C O_N\otimes\C G'\to \C O_C\otimes\C G'$ is surjective on global sections,
so that we can conclude that 
$\C S(\C O_C\otimes(\C F'-\C L),\C G)$ is the same as
$\C S(\C O_C\otimes(\C F'-\C L),\C O_C\otimes\C G')$.
\vskip\baselineskip
Now, $\C N+K_X$ is not the class of an effective divisor,
and the same will remain true if we replace $N$ by any subscheme 
of $N$ obtained by subtracting off
irreducible components of $N$.
Thus any such resulting subscheme $M$ of $N$ has the property, like $N$ itself,
that there is a component $D$ of $M$ such that 
$(M-D)\cdot D\le 1$. If $M$ is just $N$ with the reduced induced
scheme structure, then by induction on the number of components
of $M$ it follows (using \ir{keycor}) that any two components of $N$ are smooth
rational curves that are either disjoint or meet transversely
at a single point, and no sequence
$B_1$, $\ldots$, $B_i$ of distinct components exists
such that $B_i\cdot B_{1}>0$ and $B_j\cdot B_{j+1}>0$ for $1\le j<i$
(in particular, no three components meet at a single point,
and the components of $M$ form a disjoint union of trees).
\vskip\baselineskip
First assume that $N$ is reduced; i.e. that $N=N_{red}$.
Then $C$ is not a component of $N-C$.
Choose a section $\sigma_C$ of $\C O_C\otimes\C G'$,
and for each of the other components $B$ of $N$,
choose a section $\sigma_B$ of $\C O_B\otimes\C G'$
such that $\sigma_B$ does not vanish at any of the points where
$B$ meets another component of $N$. (This is possible since
$B$ is smooth and rational, so $\C O_B\otimes\C G'$
is $\C O_{\pr1}(d)$ for some $d\ge 0$, so a section can always
be chosen which does not vanish at any of a given
finite set of points of $B$.) Since $N$ has no cycles
and the components meet transversely,
it is clear that starting from $\sigma_C$
one can patch together the sections $\sigma_B$
to get a section $\sigma$ of \C G which restricts
to $\sigma_C$. Thus $\C O_N\otimes\C G'\to \C O_C\otimes\C G'$
is surjective on global sections.
\vskip\baselineskip
Now assume that $N$ is not reduced. Let $M$ be the union
of the components of $N$ which have multiplicity greater than 1 (taken with the same
multiplicities as they have in $N$) together with those multiplicity 1 components
of $N$ that meet one of these. No multiplicity 1 component
$B$ of $M$ satisfies $B\cdot (M-B)\le 1$, so there must be a
component $B$ of multiplicity more than 1 that does, and hence we also have
$B\cdot (N-B)\le 1$ for some component $B$ of $N$ of multiplicity more than 1. Now from this and
$0\to \C O_B(-N+B)\otimes \C G\to \C O_N\otimes\C G'\to \C O_J\otimes\C G'\to0$,
where $J=N-B$, we see $h^1(B, \C O_B(-N+B)\otimes \C G) = 0$, so
$\C O_N\otimes\C G'\to \C O_J\otimes\C G'$ is surjective on global sections.
But $J$ still has $C$ as a component, because either $C$ has multiplicity 1 in $N$ (and hence
$C\ne B$), or $C$ has multiplicity more than 1 in $N$ (and so even if $B=C$, $C$ remains
a component of $N-B=J$). By induction on the number of components,
we conclude that $\C O_N\otimes\C G'\to \C O_{N_{red}}\otimes\C G'$
is surjective on global sections. But $C$ is still a component of
$N_{red}$, and $\C O_{N_{red}}\otimes\C G'\to \C O_{C}\otimes\C G'$
is surjective on global sections from above, hence so is 
$\C O_{N}\otimes\C G'\to \C O_{C}\otimes\C G'$.
\smallskip}
\vskip\baselineskip

Since $\C S(\C O_L\otimes\C F,\C O_L\otimes\C G)=0$ by induction,
it suffices to show
$\C S(\C O_C\otimes(\C F'-\C L),\C O_C\otimes\C G')=0$.
If $C\cdot(\C F'-\C L)\ge 0$, then the latter is 0 (as in the previous paragraph).
Otherwise, we must have $0=\C F'\cdot C=\C G'\cdot C$ and
$C\cdot L=1$, so
$\C O_C(-1)=\C O_C\otimes(\C F'-\C L)$ and
$\C O_C=\C O_C\otimes\C G'$, which means 
$h^0(\C O_C,\C O_C\otimes(\C F'+\C G'-\C L))=0$ 
and hence again $\C S(\C O_C\otimes(\C F'-\C L),\C O_C\otimes\C G')=0$. \qed

\bye